\documentclass[journal,twoside,web]{ieeecolor}
\usepackage{generic}
\usepackage{cite}
\usepackage{amsmath,amssymb,amsfonts}
\usepackage{algorithmic}
\usepackage[dvipdfmx]{graphicx}
\usepackage{textcomp}

\usepackage{bm}
\usepackage{ascmac}
\usepackage{here}
\usepackage{color}
\usepackage{mathrsfs}
\usepackage[subrefformat=parens]{subcaption}

\usepackage{mathtools}

\newtheorem{problem}{Problem}
\newtheorem{proposition}{Proposition}
\newtheorem{definition}{Definition}

\newtheorem{theorem}{Theorem}
\newtheorem{corollary}{Corollary}
\newtheorem{rmk}{Remark}

\newcommand{\calB}{{\mathcal B}}

\newcommand{\calF}{{\mathcal F}}

\newcommand{\calN}{{\mathcal N}}

\newcommand{\calQ}{{\mathcal Q}}    
    
\newcommand{\calS}{{\mathcal S}}

\newcommand{\calX}{{\mathcal X}}

\newcommand{\bbE}{{\mathbb E}}

\newcommand{\bbP}{{\mathbb P}}    
    
\newcommand{\bbR}{{\mathbb R}}

\newcommand{\bbZ}{{\mathbb Z}}

\newcommand{\rmc}{{\rm c}}
\newcommand{\rmd}{{\rm d}}

\newcommand{\rmr}{{\rm r}}

\newcommand{\rmH}{{\rm H}}

\newcommand{\sfP}{{\sf P}}

\newcommand{\bbra}[1]{\ensuremath{[\![#1]\!]} }  

\newcommand{\kl}[2]{D_{\rm KL} \left({#1}\|{#2}\right)}
\newcommand{\ft}{N}
\newcommand{\hlf}{\frac{1}{2}}
\newcommand{\trace}{{\rm tr}}
\newcommand{\qm}{{Q}}
\newcommand{\qmk}[1]{Q_{{#1}}}  

\newcommand{\bsigma}{{\bar{\Sigma}}}

\newcommand{\what}[1]{{\widehat{#1}}}
\newcommand{\wtilde}[1]{{\widetilde{#1}}}
\newcommand{\wcheck}[1]{{\check{#1}}}

\newcommand{\dlt}[1]{}

\newcommand{\new}{{\rm new}}

\newcommand{\sfQ}{{\sf Q}}

\newcommand{\schr}{{Schr\"{o}dinger}}

\def\BibTeX{{\rm B\kern-.05em{\sc i\kern-.025em b}\kern-.08em
    T\kern-.1667em\lower.7ex\hbox{E}\kern-.125emX}}
\begin{document}
\title{Maximum Entropy Optimal Density Control of Discrete-time Linear Systems and\\ \schr~Bridges}
\author{Kaito Ito and Kenji Kashima
\thanks{This work was supported in part by JSPS KAKENHI Grant Numbers JP21J14577, JP21H04875, and by JST, ACT-X Grant Number JPMJAX2102. }
\thanks{K. Ito is with the School of Computing, Tokyo Institute of Technology, Yokohama, Japan (e-mail: ka.ito@c.titech.ac.jp).}
\thanks{K. Kashima is with the Graduate School of Informatics, Kyoto University, Kyoto, Japan (e-mail: kk@i.kyoto-u.ac.jp). }
}

\maketitle

\begin{abstract}
	We consider an entropy-regularized version of optimal density control of deterministic discrete-time linear systems.
    Entropy regularization, or a maximum entropy (MaxEnt) method for optimal control has attracted much attention especially in reinforcement learning due to its many advantages such as a natural exploration strategy.
    Despite the merits, high-entropy control policies induced by the regularization introduce probabilistic uncertainty into systems, which severely limits the applicability of MaxEnt optimal control to safety-critical systems.
		To remedy this situation, we impose a Gaussian density constraint at a specified time on the MaxEnt optimal control to directly control state uncertainty.
		Specifically, we derive the explicit form of the MaxEnt optimal density control. In addition, we also consider the case where density constraints are replaced by fixed point constraints. Then, we characterize the associated state process as a pinned process, which is a generalization of the Brownian bridge to linear systems. Finally, we reveal that the MaxEnt optimal density control gives the so-called \schr~bridge associated to a discrete-time linear system.
\end{abstract}


\section{Introduction}
\label{sec:introduction}

\IEEEPARstart{I}{n this} paper, we address the problem of steering the state of a deterministic discrete-time linear system to a desired distribution over a finite horizon while minimizing the energy cost with an entropy regularization term for control policies. Entropy-regularized optimal control is referred to as maximum entropy (MaxEnt) optimal control and has a close connection to Kullback--Leibler (KL) control~\cite{Todorov2009,Garrabe2022,Ito2022}.
There has been considerable interest in MaxEnt optimal control especially in reinforcement learning (RL)~\cite{Haarnoja2017,Haarnoja2018,Levine2018,Ho2016}.
This is because the entropy regularization offers many advantages such as performing good exploration for RL~\cite{Haarnoja2017}, robustness against disturbances~\cite{Eysenbach2021}, and equivalence between the MaxEnt optimal control problem and an inference problem~\cite{Levine2018}, to name a few.

Although the entropy regularization brings many benefits, the resulting high-entropy policy increases the uncertainty of the state of dynamics, which severely limits the applicability of MaxEnt optimal control to safety-critical systems.
Therefore, it is important to tame state uncertainty due to the stochasticity of policies.
A straightforward approach to achieve this is to impose a hard constraint in the state distribution at a specified time.
Steering the state of a dynamical system to a desired distribution without (explicit) entropy regularization has been addressed in the literature. Especially when the state distribution has a density function, this kind of problem is called a density control problem.

{\it Literature review:}
The most closely related work to the present paper is \cite{Chen2016part1}, which studies the problem of steering a Gaussian initial density of a continuous-time linear stochastic system to a final one with minimum energy. In particular, the optimal policy is derived in explicit form, and moreover, the optimal state process is shown to be a \schr~bridge associated to a linear system.
The \schr~bridge problem introduced by \schr~\cite{Schrodinger1932}, is known to be equivalent to an entropy-regularized optimal transport problem, which has received renewed attention in data sciences~\cite{Peyre2019}.
A substantial literature is devoted to reveal the relationships between stochastic optimal control and \schr~bridges in a continuous-time setting \cite{Chen2016relation,Beghi1997,Wakolbinger1989,Blaquiere1992,Dai1991,Pra1990}.
For a discrete time and space setting, the counterpart of the density control and \schr~bridge problems is explored in \cite{Chen2017robust}.

On the other hand, in \cite{Bakolas2016,Bakolas2018,Goldshtein2017}, the density control, or covariance steering problem for discrete-time linear systems is investigated.
In \cite{Bakolas2016,Bakolas2018}, the author employs a convex relaxation technique to solve a covariance steering problem with quadratic cost and a constraint on the state and the control input.
In \cite{Goldshtein2017}, the authors focus on covariance steering with only quadratic input cost and propose an efficient numerical scheme without relaxation for computing the solution.
However, unlike the continuous-time case, the obtained optimal process is {\em not} a \schr~bridge associated to the discrete-time linear system.
Indeed, in \cite{Beghi1996}, the \schr~bridge for the discrete-time linear system driven by noise acting on all components of the state vector is derived. Nevertheless, the existence and an implementable form of the optimal control that replicates the bridge are missing, and it is not possible to obtain the bridge by a linear feedback law.

Other scenarios for density control and covariance steering have been explored, including nonlinear systems, non-quadratic costs, non-Gaussian densities~\cite{Yi2020,Caluya2022,Chen2021cov,Bakolas2020,Elamvazhuthi2019}, chance constraints~\cite{Okamoto2018}, and Wasserstein terminal cost instead of a final density constraint~\cite{Halder2016,Balci2021}.
Lastly, we refer the reader to \cite{Chen2021review} for an extensive review of this area.

{\it Contributions:}
The main contributions of this paper are listed as follows:
\begin{enumerate}
    \item We analyze the MaxEnt optimal density control problem for deterministic discrete-time linear systems whose initial density and target density are both Gaussian.
		Specifically, we show the existence and uniqueness of the optimal policy and then derive its explicit form via two Lyapunov difference equations coupled through their boundary values;
    \item As a limiting case, we consider the above problem where density constraints are replaced by point constraints and show that the induced state process is a pinned process of a linear system driven by Gaussian noise. The pinned process is a generalization of the so-called Brownian bridge; see Section~\ref{sec:pinned_process} for the definition. This result is the discrete-time counterpart of \cite{Chen2016bridge};
    \item Based on the above two results, we reveal that the optimal state process of our density control problem is a \schr~bridge associated to the discrete-time linear system. This result does not require the invertibility of noise intensity supposed in \cite{Beghi1996}.
\end{enumerate}
It is worth mentioning that this study is the first one that solves a density control problem with entropy regularization, to the best of our knowledge.

{\it Organization:}
This paper is organized as follows: In Section~\ref{sec:maxent_intro}, we briefly review the MaxEnt optimal control of linear systems without density constraints.
In Section~\ref{sec:formulation}, we provide the problem formulation and derive the coupled Lyapunov equations, which are crucial ingredients in our analysis. The existence, uniqueness, and explicit form of the optimal policy are given in Section~\ref{sec:solution}. 
In Section~\ref{sec:pinned_process}, we consider the MaxEnt point-to-point steering and characterize the associated state process as a pinned process.
In Section~\ref{sec:schrodinger}, we reveal that the MaxEnt optimal density control induces the \schr~bridge associated to the linear system.
In Section~\ref{sec:example}, we present illustrative examples.
Some concluding remarks are given in Section~\ref{sec:conclusion}.

{\it Notation:}
Let $ \bbR $ denote the set of real numbers and $ \bbZ_{>0} $ denote the set of positive integers. The set of integers $ \{k,k+1,\ldots,l\} \ (k< l) $ is denoted by $ \bbra{k,l} $.
Denote by $ \calS^n $ the set of all real symmetric $ n\times n $ matrices.
For matrices $ A, B \in \calS^{n} $, we write $ A \succ B $ (resp. $A \prec B$) if $ A - B $ is positive (resp. negative) definite. Similarly, $ A \succeq B $ means that $ A - B $ is positive semidefinite.
For $ A \succ 0 $ (resp. $ A\succeq 0 $), $ A^{1/2} $ denotes the unique positive definite (resp. semidefinite) square root.
The identity matrix is denoted by $ I $, and its dimension depends on the context.
The Moore-Penrose inverse of a matrix $ A \in \bbR^{m\times n} $ is denoted by $ A^\dagger \in \bbR^{n\times m} $.
The transpose of the inverse of an invertible matrix $ A $ is denoted by $ A^{-\top} $.
The column space of a matrix $ A $ is denoted by $ {\rm ran} (A) $.
Denote the trace and determinant of a square matrix $ A $ by $ \trace (A) $ and $ \det (A) $, respectively.
The Euclidean norm is denoted by $ \|\cdot \| $.
For $ A\succ 0 $, denote $\|x\|_A := (x^\top A x)^{1/2}$.
For vectors $ x_1,\ldots,x_m \in \bbR^n $, a collective vector $ [x_1^\top \ \cdots \ x_m^\top]^\top \in \bbR^{nm} $ is denoted by $ [x_1; \ \cdots \ ; x_m] $.
Let $(\Omega, \mathscr{F}, \bbP)$ be a complete probability space and $\bbE$ be the expectation with respect to $ \bbP $.
The KL divergence (relative entropy) between probability distributions $ \bbP_x, \bbP_y $ is denoted by $ \kl{\bbP_x}{\bbP_y} $ when it is defined. We use the same notation for probability densities $ p_x,p_y $ such as $ \kl{p_x}{p_y} $.
For an $ \bbR^n $-valued random vector $ w $, $ w\sim \calN(\mu,\Sigma) $ means that $ w $ has a multivariate Gaussian distribution with mean $ \mu \in \bbR^n $ and covariance matrix $ \Sigma $. When $ \Sigma \succ 0 $, the density function of $ w\sim \calN (\mu,\Sigma) $ is denoted by $ \calN (\cdot | \mu, \Sigma) $.

\section{Brief Introduction to Maximum Entropy Linear Quadratic Optimal Control}\label{sec:maxent_intro}
We start with a brief introduction to the MaxEnt optimal control of linear systems with quadratic costs. The problem formulation is given as follows.
\begin{problem}\label{prob:LQR_general}
	Find a policy $ \pi = \{\pi_k\}_{k=0}^{\ft-1} $ that solves
	\begin{align}
				&\underset{\pi }{\rm minimize} &&  \bbE \biggl[ \frac{1}{2} (x_\ft - \bar{x}_\ft)^\top F (x_\ft - \bar{x}_\ft) \nonumber\\
				& &&\quad +  \sum_{k=0}^{\ft-1} \biggl( \frac{1}{2} \|u_k\|^2  - \varepsilon \rmH\left(\pi_k (\cdot | x_k)\right)   \biggr)  \biggr]  \label{eq:cost_general}\\
				&\text{subject to} && x_{k+1} = A_k x_k + B_k u_k, \label{eq:constraint_LQR_general}\\
				& && u_k \sim \pi_k (\cdot | x) \ \text{given} \ x_k = x, \nonumber
	\end{align}
	where $ \varepsilon > 0 ,\ x_k,\bar{x}_\ft \in \bbR^n, u_k \in \bbR^m , F \in \calS^{n} $, and $ \ft\in \bbZ_{>0} $ is the terminal time. The initial state $ x_0 $ has a finite second moment, and $ \bar{x}_\ft $ is a constant. A stochastic policy $ \pi_k (\cdot | x) $ denotes the conditional density of $ u_k $ given $ x_k = x $, and $ \rmH(\pi_k (\cdot |x)) := - \int_{\bbR^m} \pi_k (u | x) \log \pi_k (u | x) \rmd u $ denotes the entropy of $ \pi_k (\cdot | x)$.
	\hfill $ \diamondsuit $
\end{problem}
For larger values of the regularization parameter $ \varepsilon $, the entropy of the optimal policy becomes higher.

Under the invertibility of $ A_k $, define the state-transition matrix
\begin{equation}\label{eq:transition_matrix}
	\Phi (k,l) :=
	\begin{cases}
		A_{k-1} A_{k-2} \cdots A_l , & k > l ,\\
		I, & k = l  ,\\
		A_k^{-1} A_{k+1}^{-1} \cdots A_{l-1}^{-1}, &  k < l.
	\end{cases} 
\end{equation}
Similar to the conventional LQR problem~\cite{Lewis2012}, the optimal policy of Problem~\ref{prob:LQR_general} can be obtained in closed form.
The proof is shown in Appendix~\ref{app:maxent_LQR}.
\begin{proposition}\label{prop:LQR_general}
	Assume that $ \Pi_k\in \calS^{n} $ satisfies $ I + B_k^\top \Pi_{k+1} B_k \succ 0 $ for any $k\in \bbra{0,\ft-1}$ and is a solution to the following Riccati difference equation:
	\begin{align}
		 &\Pi_k = A_k^\top \Pi_{k+1} A_k - A_k^\top \Pi_{k+1} B_k (I + B_k^\top \Pi_{k+1} B_k)^{-1} \nonumber\\
		 &\hspace{2.5cm} \times B_k^\top \Pi_{k+1} A_k, ~~ k\in \bbra{0,\ft-1}, \label{eq:riccati_general}\\
		 &\Pi_\ft = F. \label{eq:riccati_terminal}
	\end{align}
	Assume further that $ A_k $ is invertible for any $ k\in \bbra{0,\ft-1} $.
	Then, the unique optimal policy of Problem~\ref{prob:LQR_general} is given by
	\begin{align}
		 \pi_{k}^* ( u | x) &= \calN \bigl(u \bigl| \ -(I + B_k^\top \Pi_{k+1} B_k)^{-1} B_k^\top \Pi_{k+1} A_k \nonumber\\
		&\qquad \times (x - \Phi(k,\ft) \bar{x}_\ft ) , \  \varepsilon(I + B_k^\top \Pi_{k+1} B_k)^{-1} \bigr), \nonumber\\
		 &\hspace{2cm} k \in \bbra{0,\ft-1}, \ x\in \bbR^n, u\in \bbR^m . \label{eq:opt_policy_general}
	 \end{align}
	 In addition, if $ F \succeq 0 $, then the solution $ \{\Pi_k\} $ to \eqref{eq:riccati_general},~\eqref{eq:riccati_terminal} satisfies $ \Pi_{k} \succeq 0 $ for any $k\in \bbra{0,\ft}$.
	 Furthermore, when $ \bar{x}_\ft = 0 $, the policy \eqref{eq:opt_policy_general} with $ \Phi(k,\ft) \bar{x}_\ft = 0 $ is still optimal without the invertibility of $ A_k $.
	 \hfill $ \diamondsuit $
\end{proposition}

The mean of the optimal policy \eqref{eq:opt_policy_general} coincides with the LQ optimal controller. In other words, the optimal policy is the LQ optimal controller perturbed by independent additive Gaussian noise with zero mean and covariance matrix $ \varepsilon (I + B_k^\top \Pi_{k+1} B_k)^{-1} $, which explicitly shows the effect of the entropy regularization.
\begin{rmk}\label{rmk:eps1}
	In Problem~\ref{prob:LQR_general}, without loss of generality, we may assume $ \varepsilon = 1 $.
To see this, let $ u_k^\varepsilon := u_k/\sqrt{\varepsilon} $ and $ \pi_k^\varepsilon (\cdot | x) $ be the conditional density of $ u_k^\varepsilon $ given $ x_k = x $. Then, 
\[
	\rmH(\pi_k^\varepsilon (\cdot | x)) = \rmH (\pi_k (\cdot | x)) - \frac{m}{2} \log \varepsilon , \ \forall x\in \bbR^n .
\]
See e.g., \cite[Section~8.6]{Cover2006}.
The cost \eqref{eq:cost_general} is rewritten as
\begin{align*}
	& \varepsilon \bbE \biggl[ \frac{1}{2\varepsilon} (x_\ft - \bar{x}_\ft)^\top F (x_\ft - \bar{x}_\ft) \nonumber\\
				&\quad +  \sum_{k=0}^{\ft-1} \biggl( \frac{1}{2} \|u_k^\varepsilon\|^2  -  \rmH\left(\pi_k^\varepsilon (\cdot | x_k)\right) - \frac{m}{2} \log \varepsilon  \biggr)  \biggr] .
\end{align*}
Therefore, by replacing $ B_k, F $ by $ \sqrt{\varepsilon} B_k, F/\varepsilon $, respectively, Problem~\ref{prob:LQR_general} reduces to the case with $ \varepsilon = 1 $.
\hfill $ \diamondsuit $
\end{rmk}

\section{Problem Formulation and Preliminary Analysis}\label{sec:formulation}
Now, we formulate the problem of steering the state of a linear system to a desired distribution.
Specifically, in this paper, we consider the following optimal control problem.
\begin{problem}[MaxEnt optimal density control problem]\label{prob:density_control}
Find a policy $ \pi = \{\pi_k\}_{k=0}^{\ft-1} $ that solves
\begin{align}
		&\underset{\pi }{\rm minimize}  && \bbE \left[ \sum_{k=0}^{\ft-1} \left( \frac{1}{2} \|u_k\|^2  - \varepsilon \rmH(\pi_k (\cdot | x_k))   \right)  \right] \label{eq:cost} \\
		&\text{subject to} && x_{k+1} = A_k x_k + B_k u_k, \label{eq:system}\\
      & &&u_k \sim \pi_k (\cdot | x) \ \text{given} \ x_k = x, \nonumber\\
      & &&x_0 \sim \calN(0,\bar{\Sigma}_0), \  x_N \sim \calN(0,\bar{\Sigma}_\ft) , \label{eq:cov_constraint}
\end{align}
where $ \varepsilon > 0, \ \bar{\Sigma}_0, \bar{\Sigma}_\ft \succ 0, \ x_k \in \bbR^n, u_k \in \bbR^m, \ N\in \bbZ_{>0} $.
\hfill $ \diamondsuit $
\end{problem}

The general case (referred to as Problem~\ref{prob:density_control}$ ' $) where the density constraints \eqref{eq:cov_constraint} are replaced by non-zero mean Gaussian distributions:
\begin{equation}\label{eq:density_constraint_general}
	x_0 \sim  \calN(\bar{\mu}_0, \bar{\Sigma}_0), \ x_\ft \sim  \calN(\bar{\mu}_\ft, \bar{\Sigma}_\ft) ,
\end{equation}
is also of interest. Later, we will see that this can be solved by combining the optimal covariance steering obtained as the solution to Problem~\ref{prob:density_control} and optimal mean steering.

To tackle Problem~\ref{prob:density_control}, let us go back to Problem~\ref{prob:LQR_general} with $ \bar{x}_\ft = 0 $ (referred to as Problem~\ref{prob:LQR_general}$ ' $).
Problem~\ref{prob:LQR_general}$ ' $ has a terminal cost instead of a constraint on the density of the final state.
Assume that $ \Pi_k\in \calS^{n} $ satisfies $ I + B_k^\top \Pi_{k+1} B_k \succ 0 $ for any $k\in \bbra{0,\ft-1}$ and is a solution to the Riccati equation \eqref{eq:riccati_general} with \eqref{eq:riccati_terminal}.
By Proposition~\ref{prop:LQR_general}, the unique optimal policy of Problem~\ref{prob:LQR_general}$ ' $ is given by
\begin{align}
	 \pi_{k}^* ( u | x) &= \calN \bigl(u \bigl| \ -(I + B_k^\top \Pi_{k+1} B_k)^{-1} B_k^\top \Pi_{k+1} A_kx,  \nonumber\\
	 &\qquad\qquad \varepsilon(I + B_k^\top \Pi_{k+1} B_k)^{-1} \bigr), \nonumber\\
	 &\qquad\qquad\quad k \in \bbra{0,\ft-1}, \ x\in \bbR^n, u\in \bbR^m . \label{eq:opt_policy}
 \end{align}
The system \eqref{eq:system} driven by the policy \eqref{eq:opt_policy} is given by
\begin{align}
	x_{k+1}^* = \bar{A}_k x_k^*  +  B_kw_k^*, \ w_k^* \sim \calN(0,\varepsilon(I + B_k^\top \Pi_{k+1} B_k)^{-1}) ,\label{eq:opt_state} 
\end{align}
where $ \{w_k^*\} $ is an independent sequence and
\begin{align*}
	\bar{A}_k &:= A_k - B_k(I + B_k^\top \Pi_{k+1} B_k)^{-1} B_k^\top \Pi_{k+1} A_k \nonumber\\
	&= (I + B_k B_k^\top \Pi_{k+1})^{-1} A_k \nonumber .
\end{align*}
By the linearity of \eqref{eq:opt_state} and the Gaussianity of $ w_k^* $, if $ x_0 $ follows a Gaussian distribution, then for any $ k \in \bbra{1,\ft} $, $ x_k^* $ has a Gaussian distribution.
Suppose that $ \Sigma_k := \bbE\left[x_k^* (x_k^*)^\top \right] $ satisfies $ \Sigma_0 = \bar{\Sigma}_0, \Sigma_\ft = \bar{\Sigma}_\ft $. Then, the policy \eqref{eq:opt_policy} is the unique optimal solution of Problem~\ref{prob:density_control}. 
In fact, for any policy satisfying \eqref{eq:cov_constraint}, the terminal cost $ \bbE[x_\ft^\top F x_\ft / 2] $ for Problem~\ref{prob:LQR_general}$ ' $ takes the same value. Thus, if the policy \eqref{eq:opt_policy} is not the unique optimal solution of Problem~\ref{prob:density_control}, it contradicts the optimality and uniqueness of \eqref{eq:opt_policy} for Problem~\ref{prob:LQR_general}$ ' $.
By \eqref{eq:opt_state}, $ \Sigma_k $ evolves as
\begin{align}\label{eq:sigma_transition}
	\Sigma_{k+1} &= \bar{A}_k \Sigma_k \bar{A}_k^\top + \varepsilon B_k(I + B_k^\top \Pi_{k+1} B_k)^{-1} B_k^\top .
\end{align}
Therefore, if $ \Sigma_0 \succ 0 $ and $ A_k $ is invertible for any $ k\in \bbra{0, \ft-1} $, it holds $ \Sigma_{k} \succ 0 $ for any $k \in \bbra{1,\ft}$. Henceforth, we assume the invertibility of $ A_k $.

Now, inspired by \cite{Chen2016part1}, we introduce $ H_k := \varepsilon\Sigma_k^{-1} - \Pi_k $.
Assume that $ \Pi_k $ and $ H_k $ are invertible on the time interval $\bbra{0,\ft}$. Noting that
\begin{align*}
	(\Pi_k + H_k)^{-1} = \Pi_k^{-1} - \Pi_k^{-1} ( \Pi_k^{-1} + H_k^{-1})^{-1} \Pi_k^{-1},
\end{align*}
we obtain
\begin{equation*}
   \Pi_k^{-1} + H_k^{-1} =\left(\Pi_k -  \frac{1}{\varepsilon} \Pi_k \Sigma_k \Pi_k \right)^{-1} .
\end{equation*}
In addition, by using \eqref{eq:riccati_general},~\eqref{eq:sigma_transition} and noting that $ A_k^\top \Pi_{k+1} \bar{A}_k = \Pi_{k} $, we get
\begin{align*}
	&A_k^{-1} (\Pi_{k+1}^{-1} + H_{k+1}^{-1}) A_k^{-\top} \nonumber \\
	&= \left( A_k^\top \left(\Pi_{k+1} - \frac{1}{\varepsilon} \Pi_{k+1} \Sigma_{k+1} \Pi_{k+1} \right) A_k		\right)^{-1} \nonumber\\
	&= \left(\Pi_k - \frac{1}{\varepsilon} \Pi_k \Sigma_k \Pi_k \right)^{-1} = \Pi_k^{-1} + H_k^{-1} .
\end{align*}
Hence, it holds
\begin{equation}\label{eq:pi_h_lyap}
	\Pi_{k+1}^{-1} + H_{k+1}^{-1} = A_k (\Pi_k^{-1} + H_k^{-1}) A_k^\top .
\end{equation}
Moreover, the Riccati equation \eqref{eq:riccati_general} is rewritten as
\begin{equation}\label{eq:pi_lyap}
	\Pi_{k+1}^{-1} = A_k\Pi_k^{-1} A_k^\top - B_k B_k^\top .
\end{equation}
By \eqref{eq:pi_h_lyap} and \eqref{eq:pi_lyap}, it holds
\begin{equation}
	H_{k+1}^{-1} = A_k H_k^{-1} A_k^\top + B_k B_k^\top .
\end{equation}
Therefore, $ P_k := H_k^{-1}$ and $Q_k := \Pi_k^{-1} $ satisfy the Lyapunov difference equations
\begin{subequations}\label{eq:lyapunov}
	\begin{align}
		P_{k+1} &= A_k P_k A_k^\top + B_k B_k^\top, \label{eq:lyapunov_P} \\
		Q_{k+1} &= A_k Q_k A_k^\top - B_k B_k^\top   \label{eq:lyapunov_Q}
	\end{align}
\end{subequations}
for $ k \in \bbra{0,\ft-1} $, and the boundary conditions $ \Sigma_0 = \bar{\Sigma}_0, \Sigma_\ft = \bar{\Sigma}_\ft  $ are written as
\begin{subequations}\label{eq:boundary}
	\begin{align}
		&\varepsilon \bar{\Sigma}_0^{-1} = P_0^{-1} + Q_0^{-1} , \label{eq:boundary_initial} \\
		&\varepsilon \bar{\Sigma}_\ft^{-1} = P_\ft^{-1} + Q_\ft^{-1}. \label{eq:boundary_end}
	\end{align}
\end{subequations}
In summary, we obtain the following proposition.
\begin{proposition}\label{prop:pre_optimal}
   Assume that for any $ k \in \bbra{0,\ft-1} $, $ A_k $ is invertible.
   Assume further that $ P_k$ and $Q_k $ satisfy the equations \eqref{eq:lyapunov_P},~\eqref{eq:lyapunov_Q} with the boundary conditions \eqref{eq:boundary_initial},~\eqref{eq:boundary_end} and are invertible on $ \bbra{0,\ft} $, and that it holds $ I + B_k^\top Q_{k+1}^{-1} B_k \succ 0 $ for any $ k\in \bbra{0,\ft-1} $.
   Then, the policy \eqref{eq:opt_policy} with $ \Pi_k = Q_k^{-1} $ is the unique optimal policy of Problem~\ref{prob:density_control}.
	\hfill $ \diamondsuit $
\end{proposition}

\begin{rmk}\label{rmk:cov_steering}
	Since Proposition~\ref{prop:LQR_general} does not require $ x_0 $ to be Gaussian, the argument in this section still holds when the density constraints~\eqref{eq:cov_constraint} are replaced by the constraints on mean and covariance:
	\begin{align}\label{eq:general_cov_constraint}
		\bbE[x_0] = 0, \ \bbE[x_0x_0^\top] = \bsigma_0,\ \bbE[x_\ft] = 0, \ \bbE[x_\ft x_\ft^\top] = \bsigma_\ft . 
 \end{align}
	That is, the optimal policy given in Proposition~\ref{prop:pre_optimal} is also optimal for Problem~\ref{prob:density_control} whose constraints~\eqref{eq:cov_constraint} are replaced by \eqref{eq:general_cov_constraint}.
	For non-zero mean constraints, see Corollary~\ref{cor:general_mean_nongauss} in Section~\ref{sec:solution}.
	\hfill $ \diamondsuit $
\end{rmk}

\section{Solution to Maximum Entropy Optimal Density Control Problem}\label{sec:solution}
In this section, we analyze the Lyapunov equations \eqref{eq:lyapunov_P},~\eqref{eq:lyapunov_Q} with the boundary conditions \eqref{eq:boundary_initial},~\eqref{eq:boundary_end}.
For the analysis, we introduce the reachability Gramian
\begin{equation}\label{eq:reachability_gramian}
	G_{\rm r} (k_1, k_0) := \sum_{k=k_0}^{k_1 -1} \Phi(k_1, k+1) B_k B_k^\top \Phi(k_1, k+1)^\top , \ k_0 < k_1,
\end{equation}
and the controllability Gramian
\begin{equation}\label{eq:controllability_gramian}
	G_{\rm c} (k_1, k_0) := \sum_{k=k_0}^{k_1 -1} \Phi(k_0, k+1) B_k B_k^\top \Phi(k_0, k+1)^\top , \ k_0 < k_1 .
\end{equation}
Note that since
\begin{equation}\label{eq:control_reachable}
	G_\rmc (k_1,k_0) = \Phi (k_0, k_1) G_\rmr (k_1, k_0) \Phi(k_0,k_1)^\top,
\end{equation}
	if $ G_\rmr (k_1,k_0) $ is invertible, $ G_\rmc (k_1,k_0) $ is also invertible.
Now, we provide the solutions to \eqref{eq:lyapunov_P},~\eqref{eq:lyapunov_Q} with \eqref{eq:boundary_initial},~\eqref{eq:boundary_end}. The proof is shown in Appendix~\ref{app:proof_lyapunov_sol}.

\begin{proposition}\label{prop:lyap_solution}
	Assume that for any $ k\in \bbra{0,\ft-1} $, $ A_k $ is invertible, and there exists $ k_\rmr \in \bbra{1,\ft}$ such that $ G_\rmr (k,0) $ is invertible for any $ k\in \bbra{k_\rmr, \ft} $ and $ G_\rmr (\ft,k) $ is invertible for any $ k\in \bbra{0,k_\rmr -1} $.
	Assume further that for
	\begin{align}
		&\hspace{-0.2cm}S_0 := \frac{1}{\varepsilon} G_{\rm c} (\ft,0)^{-\hlf} \bar{\Sigma}_0 G_{\rm c}(\ft,0)^{-\frac{1}{2}}, \\
		&\hspace{-0.2cm}S_\ft := \frac{1}{\varepsilon} G_{\rm c} (\ft,0)^{-\frac{1}{2}} \Phi(0,\ft) \bar{\Sigma}_\ft \Phi(0,\ft)^\top G_{\rm c} (\ft,0)^{-\frac{1}{2}},
	\end{align}
	the following two matrices are invertible.
	\begin{align}
		&\calF (S_0,S_\ft) := S_0 + \frac{1}{2} I - \left(S_0^\hlf S_\ft S_0^\hlf + \frac{1}{4} I\right)^\hlf,\label{eq:exception_forward}\\
		&\calB (S_0,S_\ft) :=  - S_0 + \frac{1}{2} I +\left(S_0^\hlf S_\ft S_0^\hlf + \frac{1}{4} I\right)^\hlf .\label{eq:exception_backward}
	\end{align}
   Then, the equations \eqref{eq:lyapunov_P},~\eqref{eq:lyapunov_Q} with the boundary conditions \eqref{eq:boundary_initial},~\eqref{eq:boundary_end} have two sets of solutions $ (P_{\pm,k}, Q_{\pm, k}), \ k\in \bbra{0,\ft} $ specified by
	\begin{align}
		Q_{\pm,0} &= G_{\rm c} (\ft,0)^{\frac{1}{2}} S_0^{\frac{1}{2}} \biggl( S_0 + \frac{1}{2} I \pm \left( S_0^\frac{1}{2} S_\ft S_0^\frac{1}{2} + \frac{1}{4} I  \right)^{\frac{1}{2}}  \biggr)^{-1} \nonumber\\
      &\quad \times S_0^\frac{1}{2} G_{\rm c} (\ft,0)^{\frac{1}{2}}, \label{eq:Q0} \\
		P_{\pm , 0} &= (\varepsilon\bar{\Sigma}_0^{-1} - Q_{\pm,0}^{-1} )^{-1} .
	\end{align}
	In addition, the two sets of solutions $ (P_{\pm,k}, Q_{\pm, k}) $ have the following properties.
	\begin{itemize}
		\item[(i)]  $ P_{-,k} $ and $ Q_{-,k} $ are both invertible on $ \bbra{0,\ft} $, and for any $ k\in \bbra{0,\ft-1} $, it holds $ I + B_k^\top Q_{-,k+1}^{-1} B_k \succ 0$;
		\item[(ii)] If $ Q_{+,k} $ is invertible on $ \bbra{0,\ft} $, there exists $ s \in \bbra{0,\ft-1} $ such that $ I + B_s^\top Q_{+,s+1}^{-1} B_s $ is not positive definite.
		\hfill $ \diamondsuit $
	\end{itemize}
\end{proposition}

Proposition~\ref{prop:lyap_solution} says that although a pair of Lyapunov equations \eqref{eq:lyapunov} with \eqref{eq:boundary} has two sets of solutions, only the pair $ (P_{-,k}, Q_{-,k}) $ is qualified for the construction of the optimal policy based on Proposition~\ref{prop:pre_optimal}.
In summary, we come to the main result of this section.
\begin{theorem}\label{thm:density_opt}
   Suppose that the assumptions of Proposition~\ref{prop:lyap_solution} are satisfied.
   Then, the unique optimal policy of Problem~\ref{prob:density_control} is given by
	\begin{align}
		\pi_k^{*} (u | x) &= \calN \bigl(u | - (I + B_k^\top Q_{-,k+1}^{-1} B_k)^{-1} B_k^\top Q_{-,k+1}^{-1} A_kx, \nonumber\\
      &\qquad\qquad \varepsilon(I + B_k^\top Q_{-,k+1}^{-1} B_k)^{-1} \bigr), \nonumber\\
      &\qquad \qquad\quad k\in \bbra{0,\ft-1}, \ x\in \bbR^n, u\in \bbR^m,
	\end{align}
   where $ Q_{-,k} $ is a solution to \eqref{eq:lyapunov_Q} with the initial value $ Q_{-,0} $ in \eqref{eq:Q0}.
	\hfill $ \diamondsuit $
\end{theorem}

We now give some remarks on the assumptions in Proposition~\ref{prop:lyap_solution}.
\begin{rmk}\label{rmk:inv_A}
	In many situations, $ A_k $ is expected to be invertible.
	For example, if the system~\eqref{eq:system} is obtained by a zero-order hold discretization of a continuous-time system, $ A_k $ is always invertible.
	\hfill $ \diamondsuit $
\end{rmk}

\begin{rmk}\label{rmk:inv_Gramian}
	For a time-invariant system $ A_k \equiv A, B_k \equiv B $, the reachability Gramian $ G_\rmr (k_1,k_0) $ is given by
	\[
		G_{\rm r, TI} (k_1-k_0) := \sum_{k=0}^{k_1-k_0 -1} A^k BB^\top (A^\top)^{k} .
	\]
	Then, the invertibility assumption on the reachability Gramian in Proposition~\ref{prop:lyap_solution} means that there exists $ k_\rmr \in \bbra{1,\ft} $ such that $ G_{\rm r, TI} (k) $ is invertible for any $ k\in \bbra{k_\rmr,\ft} $ and $ G_{\rm r, TI}(\ft-k) $ is invertible for any $ k\in \bbra{0,k_\rmr-1} $. In addition, it is well-known that if the system $ (A,B) $ is reachable, there exists $ \bar{k}_\rmr \in \bbra{1,n} $ such that $ G_{\rm r, TI} (k) $ is invertible for any $ k \ge \bar{k}_\rmr $. Recall that $ n $ denotes the dimension of the state space. Therefore, when $ \ft \ge \bar{k}_\rmr $ and $ \ft - \bar{k}_\rmr + 1 \ge \bar{k}_\rmr  $, that is, when $ \ft  \ge 2 \bar{k}_\rmr -1 $, the invertibility assumption on the reachability Gramian holds. In particular, taking a terminal time $ \ft \ge 2n -1 $ always ensures the invertibility.
	\hfill $ \diamondsuit $
\end{rmk}

\begin{rmk}\label{rmk:inv_FB}
	Note first that since $ \calF(S_0,S_\ft) $ and $ \calB(S_0,S_\ft) $ are continuous in the parameters $ \{A_k\}_k ,\{B_k\}_k, \bsigma_0, \bsigma_\ft $, they become singular only in exceptional cases.
	In order to consider the implication of the invertibility assumption of $ \calF(S_0,S_\ft)$ and $\calB(S_0,S_\ft) $, we consider the following two conditions: 
	\begin{align}
		\bar{\Sigma}_\ft &= \Phi (\ft,0) \bar{\Sigma}_0 \Phi (\ft,0)^\top + \varepsilon G_\rmr (\ft,0), \label{eq:forward_invertibility}\\
		\bar{\Sigma}_0 &= \Phi (0,\ft) \bar{\Sigma}_\ft \Phi (0,\ft)^\top + \varepsilon G_\rmc (\ft,0) . \label{eq:backward_invertibility}
	\end{align}
	Under \eqref{eq:forward_invertibility}, $ \calF(S_0, S_\ft) = 0 $ and under \eqref{eq:backward_invertibility}, $ \calB(S_0, S_\ft) = 0$.
	When \eqref{eq:forward_invertibility} holds, the policy \eqref{eq:opt_policy} with the terminal weight $ F = 0 $, i.e., $ \Pi_k \equiv 0, \ u_k \sim \calN(0,\varepsilon I), \forall k $, steers the state from $ \calN(0,\bar{\Sigma}_0) $ to $ \calN(0,\bar{\Sigma}_\ft) $. Hence, by the same argument as in Section~\ref{sec:formulation}, this policy is optimal for Problem~\ref{prob:density_control} although $ \Pi_k $ is not invertible.
	On the other hand, consider the time-reversed system of \eqref{eq:system}:
	\begin{align*}
		&x_{k}^- = A_k^{-1} x_{k+1}^- - A_k^{-1} B_k u_k^- , \ k = \ft-1, \ldots, 0, \\
		&x_\ft^- \sim \calN (0, \bar{\Sigma}_\ft) .
	\end{align*}
	Then, the condition \eqref{eq:backward_invertibility} indicates that an independent Gaussian noise process $ u_k^- \sim \calN(0,\varepsilon I) $ steers $ x_k^- $ from $ \calN(0,\bar{\Sigma}_\ft) $ to $ \calN(0, \bar{\Sigma}_0) $.

	When the invertibility assumption is not fulfilled, an analysis directly via the Riccati equation~\eqref{eq:riccati_general}, rather than the Lyapunov equations \eqref{eq:lyapunov}, will be appropriate as in \cite{Chen2018}. This analysis is beyond the scope of this paper.
	We would like to notice that the continuous-time counterpart of Proposition~\ref{prop:lyap_solution} \cite[Proposition~4]{Chen2016part1} also requires the invertibility of $ \calF(S_0,S_\ft)$ and $ \calB(S_0,S_\ft) $, which is not mentioned in the statement.
	\hfill $ \diamondsuit $
\end{rmk}

	Next, we extend the result of Theorem~\ref{thm:density_opt} to Problem~\ref{prob:density_control}$ ' $ with general mean distributions~\eqref{eq:density_constraint_general}.
	To this end, let us decompose $ u_k $ as $ u_k = \bar{u}_k + \wcheck{u}_k $ where $ \bar{u}_k := \bbE[u_k] $.
	Then, the evolution of the mean $ \mu_k := \bbE[x_k] $ is governed by
	\begin{equation*}
		\mu_{k+1} = A_{k} \mu_{k} + B_k \bar{u}_k, \ k\in \bbra{0,\ft-1} ,
	\end{equation*}
	which implies that $ \{\mu_k\} $ depends only on the deterministic input sequence $ \{\bar{u}_k\} $. On the other hand, the stochastic control process $ \{\wcheck{u}_k \} $ affects $ \wcheck{x}_k := x_k - \mu_k $ as
	\[
		\wcheck{x}_{k+1} = A_k \wcheck{x}_k + B_k \wcheck{u}_k , \ k\in \bbra{0,\ft-1} .
	\]
	Let $ \wcheck{\pi}_k (\cdot | \wcheck{x}) $ be the conditional density function of $ \wcheck{u}_k $ given $ \wcheck{x}_k = \wcheck{x} $. Since for any fixed $ \{\bar{u}_k\} $,
	\begin{align}
		\pi_k(\bar{u}_k + \wcheck{u} | x) &= \wcheck{\pi}_k (\wcheck{u} | x - \mu_k ) , \nonumber\\
		&\forall k \in \bbra{0,\ft-1}, \ \forall x\in \bbR^n, \forall \wcheck{u} \in \bbR^m , \label{eq:pi_picheck}
	\end{align}
	it follows that for all $ k\in \bbra{0,\ft-1} $ and $ x\in \bbR^n $,
	\begin{align}
		\rmH(\pi_k (\cdot | x)) &= - \int_{\bbR^m} \pi_k (\bar{u}_k + \wcheck{u} | x) \log \pi_k (\bar{u}_k + \wcheck{u} | x) \rmd \wcheck{u} \nonumber\\
		&= - \int_{\bbR^m} \wcheck{\pi}_k (\wcheck{u} | x - \mu_k ) \log \wcheck{\pi}_k ( \wcheck{u} | x- \mu_k ) \rmd \wcheck{u} \nonumber \\
		&= \rmH (\wcheck{\pi}_k (\cdot | x - \mu_k )) . \label{eq:entropy_check}
	\end{align}
	On the other hand, we have
	\begin{equation}\label{eq:u_ucheck}
		\bbE[ \| u_k \|^2] = \| \bar{u}_k \|^2 + \bbE[ \| \wcheck{u}_k \|^2 ] . 
	\end{equation}
	For the shifted state~$ \wcheck{x}_k $, the density constraints~\eqref{eq:density_constraint_general} are rewritten as
	\begin{equation}\label{eq:constraint_transformed}
		\wcheck{x}_0 \sim \calN(0,\bar{\Sigma}_0), \  \wcheck{x}_\ft \sim \calN(0,\bar{\Sigma}_\ft) .
	\end{equation}
	By \eqref{eq:entropy_check}--\eqref{eq:constraint_transformed}, Problem~\ref{prob:density_control}$ ' $ can be decomposed into mean steering and covariance steering as follows.
	\begin{problem}\label{prob:transformed}
		Find a deterministic control process $ \bar{u} = \{\bar{u}_k\}_{k=0}^{\ft-1} $ and a policy $ \wcheck{\pi} = \{\wcheck{\pi}_k \}_{k=0}^{\ft-1} $ that solve
		\begin{align}
			&\underset{\ \bar{u}, \wcheck{\pi} }{\rm minimize}  && \sum_{k=0}^{\ft-1} \frac{1}{2} \|\bar{u}_k \|^2 \nonumber\\
			& &&+  \bbE \left[ \sum_{k=0}^{\ft-1} \left( \frac{1}{2} \|\wcheck{u}_k\|^2  - \varepsilon \rmH(\wcheck{\pi}_k (\cdot | \wcheck{x}_k))   \right)  \right]  \label{eq:cost_transformed}\\
			&\text{subject to} && \mu_{k+1} = A_k \mu_k + B_k \bar{u}_k, \nonumber \\
			& &&\mu_0 = \bar{\mu}_0, \ \mu_\ft = \bar{\mu}_\ft, \nonumber \\
				& &&\wcheck{x}_{k+1} = A_k \wcheck{x}_k + B_k \wcheck{u}_k , \nonumber\\
				& &&\wcheck{u}_k \sim \wcheck{\pi}_k (\cdot | \wcheck{x}) \ \text{given} \ \wcheck{x}_k = \wcheck{x} ,\nonumber\\
				& &&\wcheck{x}_0 \sim \calN(0,\bar{\Sigma}_0), \  \wcheck{x}_\ft \sim \calN(0,\bar{\Sigma}_\ft) .\nonumber
		\end{align}
		\hfill $ \diamondsuit $
	\end{problem}

	Noting that $ \{\bar{u}_k\} $ is relevant only for the first term of \eqref{eq:cost_transformed}, the unique optimal solution $ \{\bar{u}_k^*\} $ to Problem~\ref{prob:transformed} is given by
	\begin{equation}\label{eq:deterministic_part}
		\bar{u}_k^* = B_k^\top \Phi (\ft,k+1)^\top G_\rmr (\ft,0)^{-1} (\bar{\mu}_\ft - \Phi (\ft,0) \bar{\mu}_0 ) .
	\end{equation}
	See also \eqref{eq:fixed_end_point} in Section~\ref{sec:pinned_process}.
	Moreover, by Theorem~\ref{thm:density_opt}, the optimal policy $ \wcheck{\pi}^* = \{\wcheck{\pi}_k^* \} $ is given by
	\begin{align*}
		\wcheck{\pi}_k^* (\wcheck{u} | \wcheck{x}) &= \calN \bigl(\wcheck{u} | - (I + B_k^\top Q_{-,k+1}^{-1} B_k)^{-1} B_k^\top Q_{-,k+1}^{-1} A_k\wcheck{x}, \nonumber\\
		&\qquad\qquad \varepsilon(I + B_k^\top Q_{-,k+1}^{-1} B_k)^{-1} \bigr) .
	\end{align*}
	Finally, by \eqref{eq:pi_picheck}, we obtain the optimal policy of Problem~\ref{prob:density_control}$ ' $ as follows.
	\begin{corollary}\label{cor:general_mean}
		Suppose that the assumptions of Proposition~\ref{prop:lyap_solution} are satisfied.
		Let $ Q_{-,k} $ be a solution to \eqref{eq:lyapunov_Q} with the initial value $ Q_{-,0} $ in \eqref{eq:Q0}.
   Then, the unique optimal policy of Problem~\ref{prob:density_control}$ ' $ is given by
		\begin{align}
			\pi_k^* (u | x) &= \calN \bigl(u | - (I + B_k^\top Q_{-,k+1}^{-1} B_k)^{-1} B_k^\top Q_{-,k+1}^{-1} A_k \nonumber\\
			&\qquad \times (x - \mu_k^* ) + \bar{u}_k^* , \ \varepsilon(I + B_k^\top Q_{-,k+1}^{-1} B_k)^{-1} \bigr), \nonumber\\
			&\qquad \qquad\quad k\in \bbra{0,\ft-1}, \ x\in \bbR^n, u\in \bbR^m, \label{eq:general_opt_policy}
		\end{align}
		where $ \{\bar{u}_k^*\} $ is given by \eqref{eq:deterministic_part} and
		\begin{equation}
			\mu_k^* := 
			\begin{cases}
			\Phi (k,0) \bar{\mu}_0 + \sum_{s=0}^{k-1} \Phi (k,s+1) B_s \bar{u}_s^* , &\\
			&\hspace{-1.5cm} k\in \bbra{1,\ft-1} , \\
			\bar{\mu}_0, &\hspace{-1.5cm}  k= 0.
			\end{cases}
		\end{equation}
		\hfill $ \diamondsuit $
	\end{corollary}

	Lastly we consider an optimal covariance steering problem, where the initial and final densities are not necessarily Gaussian. By Remark~\ref{rmk:cov_steering} and the decomposition into mean steering and covariance steering given in Problem~\ref{prob:transformed}, we obtain the following result as a straightforward consequence of Corollary~\ref{cor:general_mean}.
	\begin{corollary}\label{cor:general_mean_nongauss}
		Suppose that the assumptions of Proposition~\ref{prop:lyap_solution} are satisfied.
   Then, the unique optimal policy of Problem~\ref{prob:density_control} whose constraints \eqref{eq:cov_constraint} are replaced by
	 \begin{align}
			\bbE[x_0] &= \bar{\mu}_0, & \bbE[(x_0 - \bar{\mu}_0)(x_0 - \bar{\mu}_0)^\top] &= \bsigma_0, \\
			\bbE[x_\ft] &= \bar{\mu}_\ft, & \bbE[(x_\ft - \bar{\mu}_\ft)(x_\ft - \bar{\mu}_\ft)^\top] &= \bsigma_\ft, 
	 \end{align}
	 is given by \eqref{eq:general_opt_policy}.
		\hfill $ \diamondsuit $
	\end{corollary}

\section{Point-to-Point Steering of Discrete-time Linear Systems}\label{sec:pinned_process}
In the previous section, we discussed the density-to-density transfer over linear systems.
In this section, we investigate the point-to-point steering of linear systems with entropy regularization. This can be seen as the limiting case where the variances of the initial and target distributions go to zero.

Intuitively, the optimal policy for the point-to-point steering can be obtained by using an infinitely large terminal weight $ F $ in Problem~\ref{prob:LQR_general}.
With this in mind, let us consider Problem~\ref{prob:LQR_general} with $ F = I/\delta, \ \delta > 0 $.
In the remainder of this paper, we assume $ \varepsilon = 1 $; see Remark~\ref{rmk:eps1}.
In addition, we assume the invertibility of $ A_k $ on $ \bbra{0,\ft-1} $.
Note that if $ G_\rmc (\ft,k+1) $ is invertible, as $ \delta \downarrow  0 $, the optimal policy $ \pi_k^* (u|x) $ in \eqref{eq:opt_policy_general} converges to
\begin{align}
    &\calN \bigl( u | -(I + B_k^\top \sfQ_{k+1}^{-1} B_k)^{-1} B_k^\top \sfQ_{k+1}^{-1} A_k (x - \Phi(k,\ft) \bar{x}_\ft), \nonumber\\
		&\qquad\quad (I + B_k^\top \sfQ_{k+1}^{-1} B_k)^{-1} \bigr), \label{eq:opt_policy_invertible}
\end{align}
where
\begin{subequations}
	\begin{align}
    \sfQ_{k+1} &= A_k \sfQ_k A_k^\top - B_k B_k^\top ,  \label{eq:lyap_q_0} \\
    \sfQ_\ft &= 0  .
\end{align}
\end{subequations}
However, this is no longer true for time $ k $ when $ \sfQ_{k+1} = G_\rmc (\ft,k+1) $ is not invertible.
Instead of \eqref{eq:opt_policy_invertible}, to obtain an expression that is valid even when $ \sfQ_{k+1} $ is not invertible, recall that the mean of \eqref{eq:opt_policy_general} coincides with the form of the LQ optimal controller. In particular, by taking $ \delta \downarrow 0 $, the associated LQ optimal control from time $ k $ to the terminal time $ \ft $ given the current state $ x_k = \bar{x}_k $ converges to the solution of the following fixed end-point deterministic optimal control problem:
\begin{equation}\label{eq:fixed_end_point}
    \begin{aligned}
       &\underset{\{u_s\}_{s=k}^{\ft-1} }{\rm minimize} && \sum_{s=k}^{\ft-1} \|u_s\|^2  \\
       &\text{subject to} && x_{s+1} = A_s x_s + B_s u_s, \ s\in \bbra{k,\ft-1}, \\
       & &&x_k = \bar{x}_k, \ x_\ft = \bar{x}_\ft .
       \end{aligned}
 \end{equation}
 Since we have
 \begin{align*}
	x_\ft &= \Phi(\ft,k) x_k \\
	&\quad + [\Phi(\ft,k+1)B_k ~~ \Phi(\ft,k+2)B_{k+1} \ \cdots \ B_{\ft-1} ] \\
	&\quad \times [u_k; \ \cdots \ ;u_{\ft-1}] ,
 \end{align*}
 the optimal solution to \eqref{eq:fixed_end_point} is
 \begin{align}
	 &[u_k; \ \cdots \ ;u_{\ft-1}] \nonumber\\
	 &= [\Phi(\ft,k+1)B_k ~~ \Phi(\ft,k+2)B_{k+1} \ \cdots \ B_{\ft-1} ]^\dagger \nonumber\\
		 &\quad \times(\bar{x}_\ft - \Phi(\ft,k) \bar{x}_k) \nonumber\\
	 &= \begin{bmatrix}
		 B_k^\top \Phi (\ft,k+1)^\top \\ \vdots \\ B_{\ft-1}^\top 
	 \end{bmatrix}
	 G_\rmr (\ft,k)^\dagger (\bar{x}_\ft - \Phi(\ft,k) \bar{x}_k) , \nonumber
 \end{align}
 which immediately gives the optimal feedback control law at time $ k $ given the current state $ x_k $:
 \begin{equation*}
	 u_k = -B_k^\top \Phi (\ft,k+1)^\top G_\rmr (\ft,k)^\dagger \Phi(\ft,k)(x_k - \Phi(k,\ft)\bar{x}_\ft ) .
 \end{equation*}
In other words, as $ \delta \downarrow 0 $, the mean of the MaxEnt optimal policy \eqref{eq:opt_policy_general} converges to
\begin{equation}
	-B_k^\top \Phi (\ft,k+1)^\top G_\rmr (\ft,k)^\dagger \Phi(\ft,k)(x - \Phi(k,\ft)\bar{x}_\ft ) .
\end{equation}

On the other hand, let $ Q_{\delta,k} $ be the solution to \eqref{eq:lyapunov_Q} with $ Q_\ft = F^{-1} = \delta I $. 
Then,
\[
	Q_{\delta,k} = \delta \Phi (k,\ft) \Phi (k,\ft)^\top + G_\rmc (\ft,k) .
\]
For the covariance matrix of the MaxEnt optimal policy \eqref{eq:opt_policy_general}, we have
\begin{align*}
	&(I + B_k^\top Q_{\delta,k+1}^{-1} B_k)^{-1} = I - B_k^\top (Q_{\delta,k+1} + B_kB_k^\top)^{-1} B_k \nonumber\\
	&= I - B_k^\top \Phi (\ft,k+1)^\top ( G_\rmr (\ft,k) + \delta I )^{-1} \Phi (\ft,k+1) B_k ,
\end{align*}
where we used \eqref{eq:lyapunov_Q} and \eqref{eq:control_reachable}.
In addition, it holds
\begin{align*}
	&\lim_{\delta \downarrow 0} [\Phi(\ft,k+1)B_k \ \cdots \ B_{\ft-1} ]^\top
	(G_\rmr (\ft,k) + \delta I)^{-1} \nonumber\\
	& = [\Phi(\ft,k+1)B_k \ \cdots \ B_{\ft-1} ]^\dagger \\
	& = [\Phi(\ft,k+1)B_k \ \cdots \ B_{\ft-1} ]^\top G_\rmr (\ft,k)^\dagger,
\end{align*}
where we used the fact $ \lim_{\delta \downarrow 0} X^\top (XX^\top + \delta I)^{-1} = X^\dagger $ for any $ X\in \bbR^{m\times n} $~\cite[Chapter~3, Ex.~25]{Ben2003}.
Therefore, we obtain
\begin{align}
	&\lim_{\delta \downarrow 0} (I + B_k^\top Q_{\delta,k+1}^{-1} B_k)^{-1} \nonumber\\
	&=I - B_k^\top \Phi(\ft,k+1)^\top G_\rmr (\ft,k)^\dagger \Phi(\ft,k+1) B_k. \label{eq:cov_converge}
\end{align}
In summary, we obtain the optimal policy for $ \delta \downarrow 0 $:
\begin{align}
	&\calN \Bigl( u | -B_k^\top \Phi (\ft,k+1)^\top G_\rmr (\ft,k)^\dagger \Phi(\ft,k) \nonumber\\
	&\qquad \times (x - \Phi(k,\ft)\bar{x}_\ft ), \nonumber\\
	&\qquad  I - B_k^\top \Phi(\ft,k+1)^\top G_\rmr (\ft,k)^\dagger \Phi(\ft,k+1) B_k \Bigr) . \label{eq:opt_policy_0}
\end{align}
The state process driven by \eqref{eq:opt_policy_0} with a fixed initial state $ \bar{x}_0 \in \bbR^n $ follows
\begin{align}
	&x_{k+1}^\bullet = \what{A}_k x_k^\bullet + B_k B_k^\top \Phi (\ft,k+1)^\top G_\rmr (\ft,k)^\dagger \bar{x}_\ft + B_k w_k^\bullet, \label{eq:pin_process} \\
	&x_0^\bullet = \bar{x}_0, \ {\rm a.s.}, \label{eq:pin_initial}
\end{align}
where
\begin{align*}
	\what{A}_k := \left(I - B_k B_k^\top \Phi (\ft,k+1)^\top G_\rmr (\ft,k)^\dagger \Phi(\ft,k+1) \right)A_k 
\end{align*}
and $ \{w_k^\bullet\} $ is an independent sequence following
\begin{equation*}
	w_k^\bullet \sim \calN(0, \ I - B_k^\top \Phi(\ft,k+1)^\top G_\rmr (\ft,k)^\dagger \Phi(\ft,k+1) B_k ) .
\end{equation*}

Especially when $ k = \ft-1 $, \eqref{eq:cov_converge} is the orthogonal projection matrix onto the null space of $ B_{\ft-1} $. Hence, if $ B_{\ft-1} \neq 0 $, the covariance \eqref{eq:cov_converge} at $ k = \ft-1 $ is degenerate, and therefore the negative entropy of \eqref{eq:opt_policy_0} at $ k= \ft-1 $ and the cost \eqref{eq:cost_general} diverge to infinity.
For this reason, the control policy \eqref{eq:opt_policy_0} cannot be characterized in terms of optimality.
Nevertheless, the policy \eqref{eq:opt_policy_0} is still implementable with finite energy.
Moreover, in the following subsections, we see that the corresponding state process $ \{x_k^\bullet\} $ can be characterized as a {\em pinned process} of the associated linear system defined as follows.
\begin{definition}[\cite{Chen2016bridge}]
	Let $ \bar{x}_0, \bar{x}_\ft \in \bbR^n $ and consider the noise-driven system:
\begin{align}
    &\widetilde{x}_{k+1} = A_k \wtilde{x}_k + B_k w_k, \ w_k \sim \calN (0,I) ,\label{eq:noise_driven_point}\\
		&\wtilde{x}_0 = \bar{x}_0, \ {\rm a.s.}, \nonumber
\end{align}
where $ \{w_k\} $ is an independent sequence.
Then, a stochastic process $ \{\chi_k\}_{k=0}^{\ft} $ is said to be a pinned process associated to \eqref{eq:noise_driven_point} from $ \bar{x}_0 $ to $ \bar{x}_\ft $ if the joint distribution of $ \{\chi_k\}_{k=0}^{\ft} $ coincides with the conditional distribution of $ \{\wtilde{x}_k \}_{k=0}^\ft $ given $ \wtilde{x}_\ft = \bar{x}_\ft $.
\hfill $ \diamondsuit $
\end{definition}

\subsection{Steering from the Zero Initial State to the Zero Terminal State}
First, we deal with the case where two end points are given by $ \bar{x}_0 = 0  $ and $ \bar{x}_\ft = 0 $.
Then, by linearity of \eqref{eq:noise_driven_point}, $ \{\wtilde{x}_k\}_{k=0}^\ft $ follows a zero-mean Gaussian distribution. Therefore, it suffices to consider only the second moment of the conditional distribution of $ \{\wtilde{x}_k\} $ given $ \wtilde{x}_\ft = 0 $.
The covariance matrix of $ [\wtilde{x}_k; \wtilde{x}_s; \wtilde{x}_\ft] \ (k<s) $ is given by
\begin{equation}\label{eq:cov_three}
	\begin{bmatrix}
		\sfP_k & \sfP_k \Phi (s,k)^\top & \sfP_k \Phi (\ft,k)^\top \\
		\Phi(s,k) \sfP_k & \sfP_s & \sfP_s \Phi(\ft,s)^\top \\
		\Phi(\ft,k) \sfP_k & \Phi(\ft,s) \sfP_s & \sfP_\ft
	\end{bmatrix},
\end{equation}
where $ \sfP_k := \bbE[\wtilde{x}_k \wtilde{x}_k^\top] $ satisfies
\begin{align}\label{eq:lyap_p_0}
	&\sfP_{k+1} = A_k \sfP_k A_k^\top + B_k B_k^\top, \ \sfP_0 = 0 .
\end{align}
Taking the Schur complement of \eqref{eq:cov_three}, we obtain the covariance of the conditional distribution of $ [\wtilde{x}_k; \wtilde{x}_s] $ given $ \wtilde{x}_\ft = 0 $ as follows:
\begin{equation*}
	\begin{bmatrix}
		\wtilde{\sfP}(k,k)  &  \wtilde{\sfP}(k,s) \\
		\wtilde{\sfP}(k,s)^\top & \wtilde{\sfP}(s,s)
	\end{bmatrix} ,
\end{equation*}
where
\[
    \wtilde{\sfP}(k,s) := \sfP_k \Phi(s,k)^\top - \sfP_k \Phi(\ft,k)^\top \sfP_\ft^{-1} \Phi(\ft,s) \sfP_s .
\]
Let $ \sfP^\bullet (k,s) := \bbE[x_k^\bullet (x_s^\bullet)^\top] $. By showing $ \sfP^\bullet (k,s) = \wtilde{\sfP}(k,s) $ for $ k \le s $, we obtain the following proposition. The proof is given in Appendix~\ref{app:bridge_point}.

\begin{proposition}\label{prop:bridge_zero}
		Assume that $ \sfP_\ft = G_\rmr (\ft,0) $ is invertible, and $ A_k $ is invertible for any $ k\in \bbra{0,\ft-1} $.
    Then, the state process $ \{x_k^\bullet \}_{k=0}^{\ft} $ following \eqref{eq:pin_process},~\eqref{eq:pin_initial} with $ \bar{x}_0 = 0,  \bar{x}_\ft = 0 $ is a pinned process associated to \eqref{eq:noise_driven_point} from $ \bar{x}_0 = 0 $ to $ \bar{x}_\ft = 0 $.
		\hfill $ \diamondsuit $
\end{proposition}

\subsection{Steering Between Arbitrary Boundary Points}
Next, we consider pinned processes between the general boundary points $ \wtilde{x}_0 = \bar{x}_0, \wtilde{x}_\ft = \bar{x}_\ft $.
The analysis of the second moment is same as for Proposition~\ref{prop:bridge_zero}.
Hence, we investigate only the first moment.
Note that since
\[
	\bbE[ \wtilde{x}_k ] = \Phi (k,0) \bar{x}_0,
\]
it holds
\begin{align*}
	\ell_k &:= \bbE[\wtilde{x}_k | \wtilde{x}_\ft = \bar{x}_\ft] \nonumber\\
    &= \Phi (k,0) \bar{x}_0 + \sfP_k \Phi (\ft,k)^\top \sfP_\ft^{-1} (\bar{x}_\ft - \Phi (\ft,0) \bar{x}_0) .
\end{align*}
From \eqref{eq:pin_process}, by showing
	\begin{equation}\label{eq:first_moment}
		\ell_{k+1} = \what{A}_k \ell_k + B_k B_k^\top \Phi (\ft,k+1)^\top G_\rmr (\ft,k)^\dagger \bar{x}_\ft,
	\end{equation}
we obtain the main result of this section. The proof is given in Appendix~\ref{app:bridge_point}.
\begin{theorem}\label{thm:point_bridge}
	Assume that $ G_\rmr (\ft,0) $ is invertible, and $ A_k $ is invertible for any $ k\in \bbra{0,\ft-1} $.
	Then, for any $ \bar{x}_0, \bar{x}_\ft \in \bbR^n $, the state process $ \{x_k^\bullet \}_{k=0}^{\ft} $ following \eqref{eq:pin_process},~\eqref{eq:pin_initial} is a pinned process associated to \eqref{eq:noise_driven_point} from $ \bar{x}_0  $ to $ \bar{x}_\ft $.
	\hfill $ \diamondsuit $
\end{theorem}

\section{\schr~Bridges of Discrete-time Linear Systems}\label{sec:schrodinger}
In this section, we explore the connection between the MaxEnt optimal density control and \schr~bridges.
Given a reference stochastic process with finite horizon $ \ft $ and two marginals, the \schr~bridge problem seeks to find the closest process to the reference in terms of the KL divergence such that its marginals at times $ k = 0 $ and $ k=\ft $ coincide with the given marginals; see e.g., \cite{Leonard2014,Chen2021liaison} for reviews.
Specifically, in this paper, the reference process and marginals are given by the noise-driven linear system \eqref{eq:noise_driven_point} and the Gaussian distributions \eqref{eq:cov_constraint}, respectively. Then, the \schr~bridge problem is formulated as follows.
\begin{problem}[\schr~bridge problem]\label{prob:schrodinger}
		Let $ \{\wtilde{x}_k\}_{k=0}^\ft $ be the solution to the linear system \eqref{eq:noise_driven_point} with the initial state $ \wtilde{x}_0 \sim \calN(0,\bar{\Sigma}_0) ,  \bsigma_0 \succ 0 $
    and denote by $ \wtilde{\bbP}_{0:\ft} $ the probability distribution of $ \{\wtilde{x}_k\}_{k=0}^\ft $ on $ \calX := (\bbR^n)^{\ft+1} $. Then, find a probability distribution $ \what{\bbP}_{0:\ft} $ on $ \calX $ that solves
	\begin{equation}\label{eq:schrodinger}
	\begin{aligned}
		\underset{\what{\bbP}_{0:\ft} \in \Pi (\bar{\Sigma}_0, \bar{\Sigma}_\ft)}{\rm minimize} ~~ &\kl{\what{\bbP}_{0:\ft}}{\wtilde{\bbP}_{0:\ft}} ,
		\end{aligned}
	\end{equation}
    where $ \what{\bbP}_{0:\ft} \in \Pi (\bar{\Sigma}_0, \bar{\Sigma}_\ft) $ means that the marginals of $ \what{\bbP}_{0:\ft} $ at $ k = 0 $ and $ k = \ft $ are equal to $ \calN(0,\bar{\Sigma}_0) $ and $ \calN(0,\bar{\Sigma}_\ft), \ \bsigma_\ft \succ 0 $, respectively.
    \hfill $ \diamondsuit $
\end{problem}

The optimal solution $ \what{\bbP}_{0:\ft} $ to the above problem is referred to as the \schr~bridge associated to \eqref{eq:noise_driven_point} from $ \calN(0,\bar{\Sigma}_0) $ to $ \calN(0,\bar{\Sigma}_\ft) $. In this paper, we also refer to the corresponding stochastic process $ \{\what{x}_k\}_{k=0}^\ft \sim \what{\bbP}_{0:\ft} $ as the \schr~bridge.
Since the KL divergence $ D_{\rm KL} (\cdot \| \wtilde{\bbP}_{0:\ft} ) $ is strictly convex, and the constraint set $ \Pi (\bar{\Sigma}_0, \bar{\Sigma}_\ft) $ is a convex subset of the vector space of all bounded measures on $ \calX $, the \schr~bridge problem is a strictly convex problem. Hence, Problem~\ref{prob:schrodinger} admits at most one solution.

Denote by $ \what{\bbP}_{0,\ft} $ and $ \what{\bbP}_{1:\ft-1|0,\ft} $ the joint distribution of $ (\what{x}_0, \what{x}_\ft) $ and the conditional distribution of $ \{\what{x}_k \}_{k=1}^{\ft-1} $ given $ ( \what{x}_0, \what{x}_\ft )$, respectively. The same notation is used for $ \wtilde{\bbP} $. For the analysis of Problem~\ref{prob:schrodinger}, we employ the decomposition of the KL divergence~\cite[Theorem~2]{Leonard2014some}:
\begin{align}
	&\kl{\what{\bbP}_{0:\ft}}{\wtilde{\bbP}_{0:\ft}} = \kl{\what{\bbP}_{0,\ft}}{\wtilde{\bbP}_{0,\ft}} \nonumber \\
		&+ \int D_{\rm KL} \Bigl(\what{\bbP}_{1:\ft-1|0,\ft} (\cdot | \what{x}_0 = \bar{x}_0, \what{x}_\ft = \bar{x}_\ft)  \| \nonumber\\
		&\hspace{0.5cm} \wtilde{\bbP}_{1:\ft-1|0,\ft} (\cdot | \wtilde{x}_0 = \bar{x}_0, \wtilde{x}_\ft = \bar{x}_\ft) \Bigr)  \what{\bbP}_{0,\ft} (\rmd \bar{x}_0 \rmd \bar{x}_\ft) .\label{eq:separate_KL}
\end{align}
The first term depends only on the joint distribution of $ (\what{x}_0, \what{x}_\ft) $.
On the other hand, by taking $ \what{\bbP}_{1:\ft-1|0,\ft} (\cdot | \what{x}_0 = \bar{x}_0, \what{x}_\ft = \bar{x}_\ft) = \wtilde{\bbP}_{1:\ft-1|0,\ft} (\cdot | \wtilde{x}_0 = \bar{x}_0, \wtilde{x}_\ft = \bar{x}_\ft) $, the second term attains its minimum value $ 0 $.
This implies that for any $ \bar{x}_0, \bar{x}_\ft\in \bbR^n $, the \schr~bridge $ \{\what{x}_k\} $ given $ \what{x}_0 = \bar{x}_0, \what{x}_\ft = \bar{x}_\ft $ is a pinned process associated to \eqref{eq:noise_driven_point} from $ \wtilde{x}_0 = \bar{x}_0 $ to $\wtilde{x}_\ft = \bar{x}_\ft $.

In \cite{Beghi1996}, the optimal solution to Problem~\ref{prob:schrodinger} is derived only for the restricted case where $ B_k $ has full row rank.
In the following theorem, for general $ B_k $, we reveal that the MaxEnt optimal state process $ \{x_k^*\}_{k=0}^\ft $ of Problem~\ref{prob:density_control} is actually an optimal solution to Problem~\ref{prob:schrodinger}.
The proof is given in Appendix~\ref{app:schrodinger_proof}.

\begin{theorem}\label{thm:schrodinger}
	Suppose that the assumptions of Proposition~\ref{prop:lyap_solution} are satisfied.
	Consider the optimal state process $ \{x_k^*\}_{k=0}^\ft $ of Problem~\ref{prob:density_control} with $ \varepsilon = 1 $ following
	\begin{align}
	&x_{k+1}^* = A_{Q,k} x_k^* + B_{Q,k} w_k, \ w_k \sim \calN(0,I), \label{eq:opt_transition}
	\end{align}
	where
	\begin{align*}
		&A_{Q,k} := A_k - B_k(I + B_k^\top Q_{-,k+1}^{-1} B_k)^{-1} B_k^\top Q_{-,k+1}^{-1} A_k, \\
		&B_{Q,k} := B_k (I + B_k^\top Q_{-,k+1}^{-1} B_k)^{-1/2} .
	\end{align*}
    Then, the probability distribution $ \bbP_{0:\ft}^* $ of $ \{x_k^*\}_{k=0}^\ft $ is the unique optimal solution to Problem~\ref{prob:schrodinger}.
		\hfill $ \diamondsuit $
\end{theorem}

\section{Illustrative Examples}\label{sec:example}
\begin{figure}[t]
	\begin{minipage}[b]{1.0\linewidth}
		\centering
		\includegraphics[keepaspectratio, scale=0.35]
		{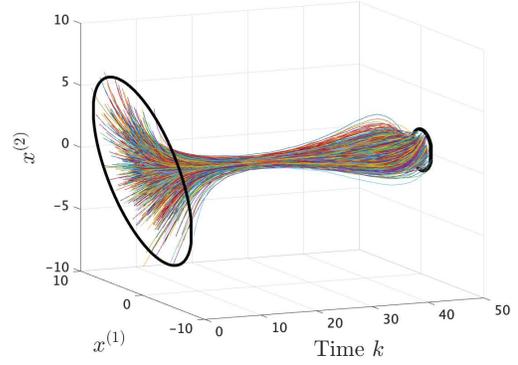}
		\subcaption{$ \varepsilon=0.02 $}\label{fig:opt_state_002}
	\end{minipage}
	\begin{minipage}[b]{1.0\linewidth}
		\centering
		\includegraphics[keepaspectratio, scale=0.35]
		{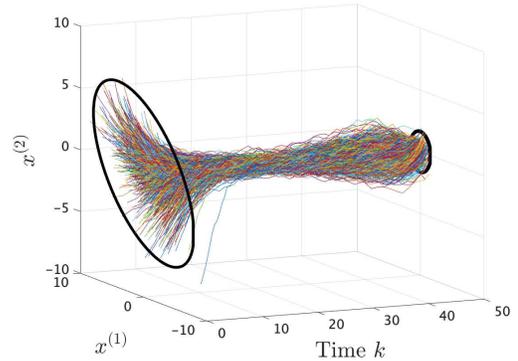}
		\subcaption{$ \varepsilon= 1.0 $}\label{fig:opt_state_1}
	\end{minipage}

	\vspace{0.5cm}
	\begin{minipage}[b]{0.49\linewidth}
		\centering
		\includegraphics[keepaspectratio, scale=0.23]
		{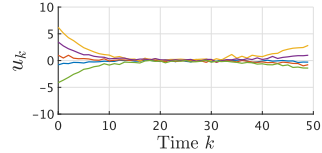}
		\subcaption{$ \varepsilon=0.02 $}\label{fig:opt_control_002}
	\end{minipage}
	\begin{minipage}[b]{0.48\linewidth}
		\centering
		\includegraphics[keepaspectratio, scale=0.23]
		{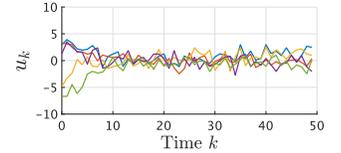}
		\subcaption{$ \varepsilon=1.0 $}\label{fig:opt_control_1}
	\end{minipage}
	\caption{$ 1000 $ samples of the optimal state process $ x_k^* = [x^{(1)}_k \ x^{(2)}_k ]^\top $ and $ 5 $ samples of the optimal control process for \eqref{eq:ex_system},~\eqref{eq:ex_cov} (colored lines) with $ \varepsilon = 0.02, 1.0 $. The black ellipses are given by \eqref{eq:ellipse}.}\label{fig:opt_trajectory}
\end{figure}

In this section, two examples illustrate the obtained results.
The system is given by
\begin{align}\label{eq:ex_system}
	A_k = 
	\begin{bmatrix}
		0.9 & 0.1\\
		0.05 & 1.2
	\end{bmatrix}, \
	B_k = [0 ~~ 0.22]^\top , \ \forall k .
\end{align}
The initial density $ \calN(0,\bar{\Sigma}_0) $ and the target density $ \calN(0,\bar{\Sigma}_\ft) $ are determined by
\begin{align}\label{eq:ex_cov}
	\bar{\Sigma}_0 =
	\begin{bmatrix}
		7 & 3\\
		3 & 5
	\end{bmatrix}, \ 
	\bar{\Sigma}_\ft = 0.3 I ,
\end{align}
and set $ N = 50 $. Then the assumptions of Proposition~\ref{prop:lyap_solution} and Theorem~\ref{thm:point_bridge} are fulfilled.
In Fig.~\ref{fig:opt_trajectory}, the sample paths of the optimal state process and the optimal control process obtained in Theorem~\ref{thm:density_opt} are shown for different $ \varepsilon $. The colored lines indicate the samples, and the black ellipses in Figs.~\ref{fig:opt_state_002},~\ref{fig:opt_state_1} are the $ 3\sigma $ covariance ellipses for $ \calN(0,\bar{\Sigma}_0) $ and $ \calN(0,\bar{\Sigma}_\ft) $ given by
\begin{equation}\label{eq:ellipse}
	\left\{ x\in \bbR^2 : x^\top \bar{\Sigma}_k^{-1} x = 3^2 \right\}, \  k= 0, \ft .
\end{equation}
As can be seen, the optimal control for $ \varepsilon = 1.0 $ has a larger variance than for $ \varepsilon = 0.02 $ resulting in higher entropy of the policy. The same is true for the state processes. Despite the difference in entropy, the state is steered to the target distribution in both cases.
Note that for $ \varepsilon = 1.0 $, the eigenvalues of $ Q_{-,\ft}^{-1} $, which corresponds to the terminal weight matrix $ F $ in Problem~\ref{prob:LQR_general}, are given by $ \{-45.81, 3.33\} $.
This clarifies that it is not enough to consider only positive semidefinite terminal weight matrices for constructing the MaxEnt optimal density control.

\begin{figure}[t]
	\begin{minipage}[b]{0.49\linewidth}
		\centering
		\includegraphics[keepaspectratio, scale=0.23]
		{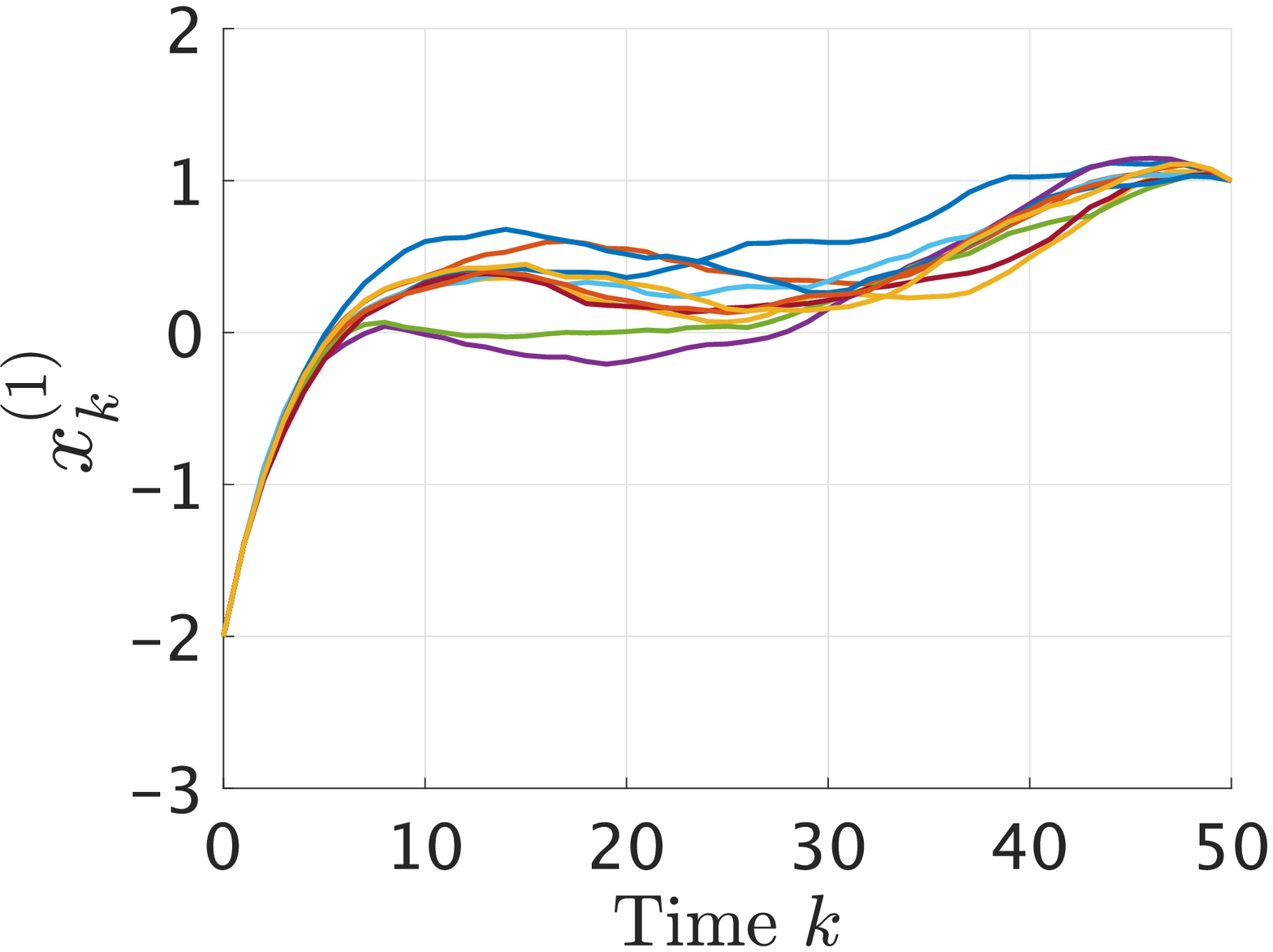}
		\subcaption{State process $ x_k^{(1)} $}\label{fig:pin_state1}
	\end{minipage}
	\begin{minipage}[b]{0.48\linewidth}
		\centering
		\includegraphics[keepaspectratio, scale=0.23]
		{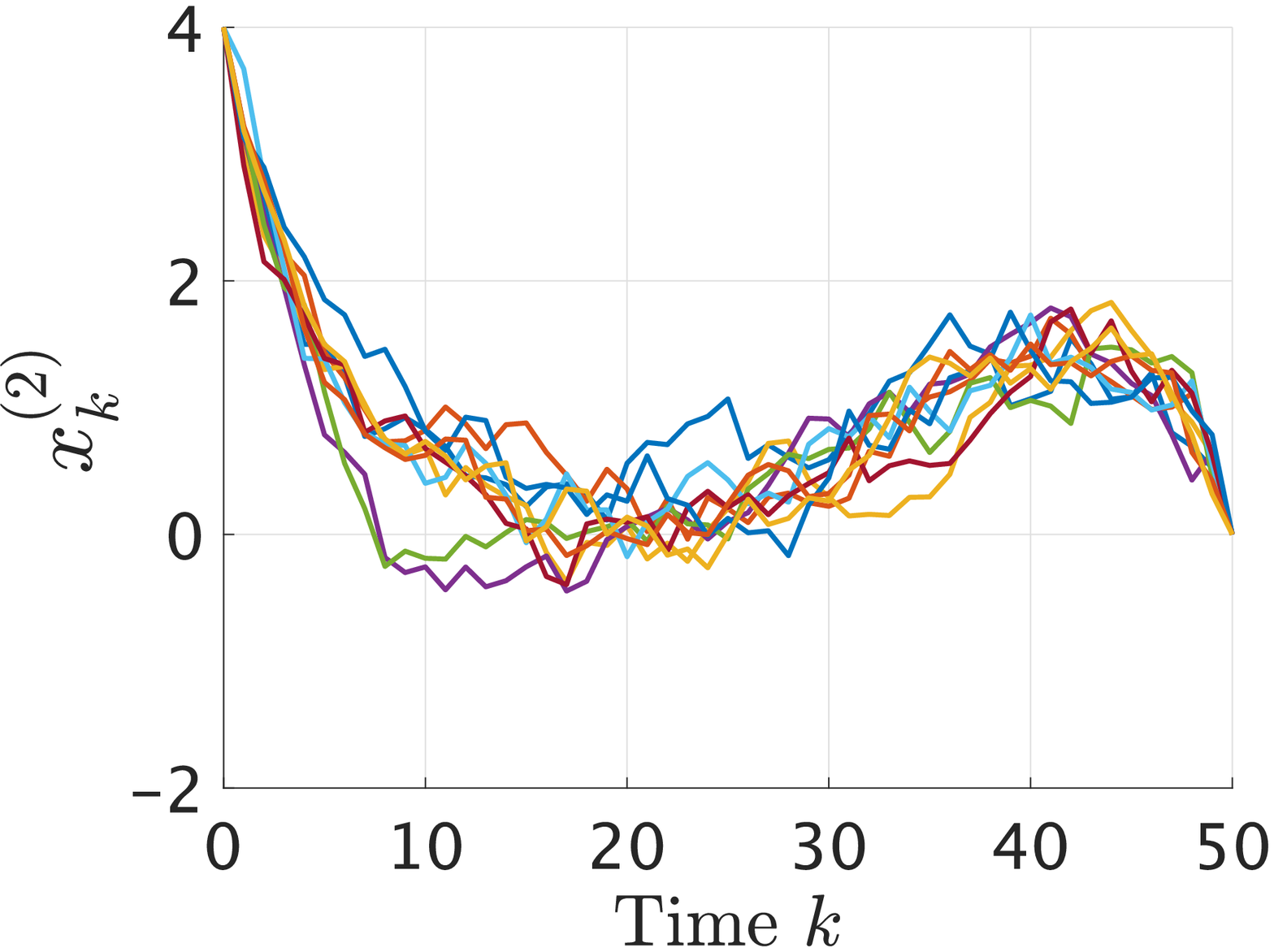}
		\subcaption{State process $ x_k^{(2)} $}\label{fig:pin_state2}
	\end{minipage}

	\vspace{0.5cm}
	\begin{minipage}[b]{1.0\linewidth}
		\centering
		\includegraphics[keepaspectratio, scale=0.26]
		{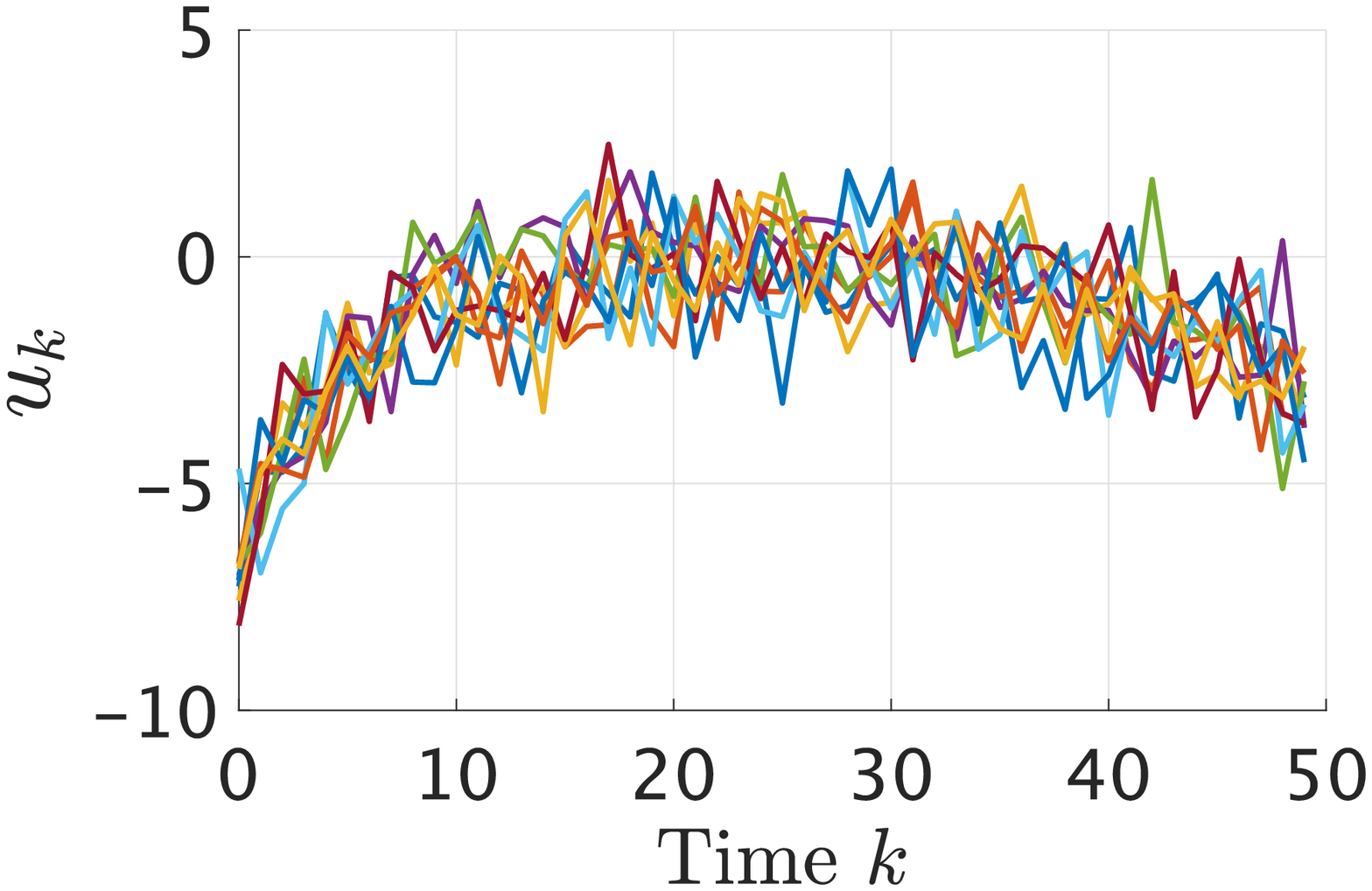}
		\subcaption{Control process $ u_k $}\label{fig:pin_u}
	\end{minipage}

	\caption{$ 10 $ samples of the pinned process $ x_k^\bullet = [x^{(1)}_k \ x^{(2)}_k ]^\top $ following \eqref{eq:pin_process} and the control process given by \eqref{eq:opt_policy_0} for \eqref{eq:ex_system}.}\label{fig:pin}
\end{figure}

Next, we describe the point-to-point steering for \eqref{eq:ex_system} with $ \bar{x}_0 = [-2 \ 4]^\top, \ \bar{x}_\ft = [1 \ 0]^\top $.
Fig.~\ref{fig:pin} depicts the samples of the pinned process $ \{x_k^\bullet\} $ following \eqref{eq:pin_process} and the corresponding control process given by \eqref{eq:opt_policy_0}.
It can be seen that all the samples of the state are transferred to the target state $ \bar{x}_\ft $ in spite of the fluctuations along the paths.

\section{Conclusion}\label{sec:conclusion}
In this paper, we analyzed MaxEnt optimal density control of deterministic discrete-time linear systems.
First, we revealed the existence and uniqueness of the optimal solution, and provided the explicit construction of the optimal policy. 
Next, we considered the point-to-point steering and characterized the associated state process as a pinned process of the noise-driven linear system.
Finally, we showed that the MaxEnt optimal density control gives the \schr~bridge associated to the discrete-time linear system.
Future work includes removing or relaxing the invertibility assumption on $ A_k, \calF(S_0,S_\ft), \calB(S_0,S_\ft) $, and the reachability Gramian in Proposition~\ref{prop:lyap_solution}. We expect that the invertibility condition on the reachability Gramian can be relaxed to the assumption that $ G_\rmr (\ft,0) $ is invertible used in Theorem~\ref{thm:point_bridge}. Another direction of future work is to extend the results to the case where a quadratic state cost is also present.
These issues may be addressed by directly analyzing the Riccati equation rather than the Lyapunov equations.

MaxEnt optimal density control of {\em stochastic} systems is also of interest.
Note that for a continuous-time stochastic system, the density control problem can be analyzed via coupled Lyapunov equations only when the input and noise channels of the system coincide~\cite{Chen2016part1}.
When the channels are distinct, the approach based on the Lyapunov equations can no longer be employed for the density control, and an alternative approach e.g., an optimization based approach~\cite{Chen2016part2}, is needed.
Similar to the continuous-time case, in the MaxEnt optimal density control problem for a deterministic discrete-time system, the input and the noise induced by entropy regularization are injected into the system through the same channel. Therefore, an optimization based approach~\cite{Bakolas2016,Bakolas2018}, rather than the one based on Lyapunov equations, would work for MaxEnt optimal density control of stochastic systems.

\appendices

\section{Proof of Proposition~\ref{prop:LQR_general}}\label{app:maxent_LQR}
Similar to the conventional optimal control, the dynamic programming approach can be applied to MaxEnt optimal control \cite{Levine2018}. Let us introduce the value function associated to Problem~\ref{prob:LQR_general}:
\begin{align*}
	&V(k,x) := \inf_{ \{\pi_s\}_{s=k}^{\ft-1} } \bbE \biggl[ \frac{1}{2} (x_\ft - \bar{x}_\ft)^\top F (x_\ft - \bar{x}_\ft)   \nonumber\\
	&+  \sum_{s=k}^{\ft-1} \biggl( \frac{1}{2} \|u_s\|^2  - \varepsilon \rmH(\pi_s (\cdot | x_s))   \biggr) \biggr| \ x_k = x \biggr], \ \text{subj. to \eqref{eq:constraint_LQR_general}} , \nonumber\\
	&\hspace{5.0cm} k\in \bbra{0,\ft-1}, \ x\in \bbR^n , \\
	&V(\ft,x) := \frac{1}{2} (x - \bar{x}_\ft)^\top F (x - \bar{x}_\ft) , ~~ x\in \bbR^n .
\end{align*}
It is known that $ V $ satisfies the {\em soft} Bellman equation:
\begin{align}
	&V(k,x) = -\varepsilon \log \left[ \int_{\bbR^m} \exp \left( - \frac{1}{\varepsilon} \calQ_k (x,u)  \right) \rmd u   \right], \nonumber\\
	&\hspace{4cm} k\in \bbra{0,\ft-1},\ x\in \bbR^n, \label{eq:bellman}
\end{align}
where the Q-function is defined by
\[
	\calQ_k (x,u) := \frac{1}{2} \| u\|^2 + V(k+1, A_k x + B_k u) , \ x\in \bbR^n,  u\in \bbR^m.
\]
Then the unique optimal policy is given by
\begin{equation}\label{eq:maxent_opt}
	\pi_k^* (u | x) \propto \exp \left(  - \frac{1}{\varepsilon} \calQ_k (x,u) \right), \ x\in \bbR^n,  u\in \bbR^m.
\end{equation}
In what follows, we calculate the Q-function analytically.
For conciseness, we drop the subscripts of $ A_k, B_k $.
First, for $ k = \ft-1 $, we have
\begin{align*}
	&\calQ_{\ft-1} (x,u) \\
	&= \frac{1}{2} \| u + (I + B^\top \Pi_\ft B)^{-1} B^\top \Pi_\ft (A x - \bar{x}_\ft) \|_{I + B^\top \Pi_\ft B}^2 \nonumber \\
	&\quad + \frac{1}{2} (x - A^{-1} \bar{x}_\ft)^\top \Pi_{\ft-1} (x - A^{-1} \bar{x}_\ft) ,
\end{align*}
where we used the invertibility of $ I + B^\top \Pi_{\ft} B $ by assumption. By combining this with \eqref{eq:maxent_opt}, we obtain the policy \eqref{eq:opt_policy_general} for $ k = \ft-1 $.
On the other hand, note that under $ I + B^\top \Pi_{\ft} B \succ 0 $, for any $ \mu \in \bbR^m $,
\begin{align*}
	&\int_{\bbR^m} \exp \left( - \frac{1}{2} \| u - \mu \|_{(I + B^\top \Pi_\ft B)/\varepsilon }^2  \right) \rmd u  \\
	&= \sqrt{ (2\pi)^m \det (\varepsilon (I  +B^\top \Pi_\ft B)^{-1} ) } .
\end{align*}
Then by \eqref{eq:bellman}, a straightforward calculation yields
\begin{align*}
	&V(\ft-1,x) \\
	&= \frac{1}{2} (x - \Phi (\ft-1,\ft) \bar{x}_\ft)^\top \Pi_{\ft-1} (x - \Phi (\ft-1,\ft) \bar{x}_\ft) \nonumber \\
	&\quad - \varepsilon \log \sqrt{ (2\pi)^m \det (\varepsilon (I  +B^\top \Pi_\ft B)^{-1} ) },
\end{align*}
where the first term takes the same form as $ V(\ft,x) $ and the second term does not depend on $ x $. 
In addition, if $ F \succeq 0 $, then $ \Pi_\ft $ has the square root $ \Pi_\ft^{1/2} $, and by \eqref{eq:riccati_general}, $ \Pi_{\ft-1} $ can be written as
\begin{equation*}
	\Pi_{\ft-1} = A^\top \Pi_{\ft}^{1/2} ( I + \Pi_{\ft}^{1/2}  B  B^\top \Pi_{\ft}^{1/2})^{-1} \Pi_{\ft}^{1/2} A \succeq 0 .
\end{equation*}
Moreover, it is obvious that the above argument holds when $ \bar{x}_\ft = 0 $ and $ A_{\ft-1} $ is not invertible.
Finally, by applying the same procedure as above for $ k = \ft-2, \ldots, 0 $, we obtain the desired result.

\section{Proof of Proposition~\ref{prop:lyap_solution}}\label{app:proof_lyapunov_sol}
First, introduce the change of variables $ \xi_k := G_{\rm c} (\ft,0)^{-1/2} \Phi(0,k) x_k $. Then the system~\eqref{eq:system} is transformed into
	\begin{equation}\label{eq:system_xi}
		\xi_{k+1}  = \xi_k + \underbrace{G_{\rm c}(\ft, 0)^{-1/2} \Phi(0,k+1) B_k}_{\eqqcolon B_{\new, k}} u_k .
   \end{equation}
   We will prove the statement in this new set of coordinates and then turn back to the original set of coordinates at the end.
   The Lyapunov equations associated with the transformed system \eqref{eq:system_xi} are given by
	 \begin{subequations}
		\begin{align}
			P_{\new, k+1} &= P_{\new, k} + B_{\new, k} B_{\new, k}^\top, \\
			Q_{\new, k+1} &= Q_{\new, k} - B_{\new, k} B_{\new, k}^\top. \label{eq:lyapunov_Qnew}
		\end{align}
	 \end{subequations}
	The relationship between $ Q_{\new,k}$ and $Q_k $ is as follows.
\begin{equation}\label{eq:Q_Qnew_relation}
	Q_{\new, k} = G_{\rm c} (\ft,0)^{-\hlf} \Phi(0,k) Q_k \Phi(0,k)^\top G_{\rm c} (\ft,0)^{-\hlf} .
\end{equation}
	Indeed, substituting \eqref{eq:Q_Qnew_relation} into $ Q_{\new, k+1} - Q_{\new,k} $ yields
\begin{align}
	&Q_{\new, k+1} - Q_{\new, k}
	=G_{\rm c}(\ft,0)^{-\hlf} \Phi(0,k) \nonumber\\
   &\quad \times(A_k^{-1} Q_{k+1} A_k^{-\top} - Q_k) \Phi(0,k)^\top G_{\rm c}(\ft,0)^{-\hlf} \nonumber\\
	&= - G_{\rm c}(\ft,0)^{-\hlf} \Phi(0,k+1) B_k B_k^\top \Phi(0,k+1)^\top G_{\rm c}(\ft,0)^{-\hlf} \nonumber\\
	&= - B_{\new, k} B_{\new, k}^\top, \nonumber
\end{align}
which coincides with \eqref{eq:lyapunov_Qnew}. 
The controllability Gramian and the reachability Gramian corresponding to $ \xi_k $ are given by
\begin{align*}
   G_{\rmc, \new} (k_1,k_0) = G_{\rmr, \new} (k_1,k_0) :=  \sum_{k = k_0}^{k_1-1} B_{\new, k} B_{\new,k}^\top
\end{align*}
satisfying $G_{{\rm c},\new} (\ft, 0) = G_{{\rm r},\new} (\ft, 0)  = I $.
   The initial and final covariance matrices for $ \xi_0 $ and $ \xi_\ft $ are given by
	\begin{align*}
		\bar{\Sigma}_{\new,0} &:=  G_{\rm c} (\ft,0)^{-1/2} \bar{\Sigma}_0 G_{\rm c} (\ft,0)^{-1/2} = \varepsilon S_0 , \\
		\bar{\Sigma}_{\new,\ft} &:=  G_{\rm c} (\ft,0)^{-1/2} \Phi(0,\ft) \bar{\Sigma}_\ft \Phi(0,\ft)^\top G_{\rm c} (\ft,0)^{-1/2} \nonumber\\
      & = \varepsilon S_\ft. 
	\end{align*}

In what follows, to simplify notation, we omit the subscript ``$ \new $.''
By substituting
\begin{equation}\label{eq:lyap_0N}
	P_{\ft} = P_{ 0} + I, \ Q_{ \ft} = Q_{ 0} - I
\end{equation}
into the boundary conditions \eqref{eq:boundary_initial},~\eqref{eq:boundary_end}, similar to the proof of \cite[Proposition~4]{Chen2016part1}, we obtain the two sets of initial values
\begin{align*}
	&Q_{\pm, 0} = S_0^{\hlf} \biggl( S_0 +\frac{1}{2} I \pm \left( S_0^{\hlf} S_\ft S_0^{\hlf} + \frac{1}{4}I  \right)^{\hlf}  \biggr)^{-1} S_0^{\hlf}, \\
   &P_{\pm, 0} = (S_0^{-1} - Q_{\pm,0}^{-1})^{-1} .
\end{align*}
Here, the invertibility of $ \calF(S_0,S_\ft) $ and $ \calB(S_0,S_\ft) $ guarantees the invertibility of $ Q_{-,0} $ and $ P_{-,0} $.
Next, we show that $ Q_{-,k} = Q_{-,0} - G_{\rmr} (k,0) $ is invertible on $ \bbra{0,\ft} $.
Note that formally
\begin{align}
   &(Q_{-,0} - G_\rmr (k,0))^{-1} = - G_\rmr (k,0)^{-1} \nonumber\\
   &\quad - G_\rmr (k,0)^{-1} \left(Q_{-,0}^{-1} - G_\rmr (k,0)^{-1}\right)^{-1} G_\rmr (k,0)^{-1} \nonumber\\
	 &= -G_\rmr (k,0)^{-1} - G_\rmr (k,0)^{-1} S_0^{\hlf} \biggl[ S_0 + \frac{1}{2} I \nonumber\\
    & - \left( S_0^{\hlf} S_\ft S_0^{\hlf} + \frac{1}{4} I  \right)^{\hlf} - S_0^{\hlf} G_\rmr (k,0)^{-1} S_0^{\hlf}  \biggr]^{-1} S_0^{\hlf} G_\rmr (k,0)^{-1} . \label{eq:Q_inv_forward}
\end{align}
By assumption, for any $ k \in \bbra{k_\rmr ,\ft} $, $ G_{\rmr }(k,0) $ is invertible. The term in the square brackets obviously attains its maximum at $ k = \ft $
\[
	\frac{1}{2} I - \left( S_0^{1/2} S_\ft S_0^{1/2} + \frac{1}{4} I  \right)^{1/2} \prec 0 .
\]
Therefore the term in the brackets is invertible for any $ k \in \bbra{k_\rmr, \ft} $. This implies that $ Q_{-,k} $ has the inverse matrix \eqref{eq:Q_inv_forward}.
On the other hand, $ Q_{-,k} $ also admits the expression $ Q_{-,k} = Q_{-,\ft} + G_{\rmr } (\ft,k) $.
Hence, the inverse matrix is formally
\begin{align}
   &\hspace{-0.3cm}(Q_{-,\ft} + G_\rmr (\ft,k))^{-1} = G_\rmr (\ft,k)^{-1} \nonumber\\
   &- G_\rmr (\ft,k)^{-1} \left(Q_{-,\ft}^{-1} + G_\rmr (\ft,k)^{-1}\right)^{-1} G_\rmr (\ft,k)^{-1} \nonumber\\
	&\hspace{-0.3cm}=G_\rmr (\ft,k)^{-1} - G_\rmr (\ft,k)^{-1} \nonumber\\
	 &\quad \times \left(( Q_{-,0} - I )^{-1} + G_\rmr (\ft,k)^{-1}\right)^{-1} G_\rmr (\ft,k)^{-1} \nonumber\\
	&\hspace{-0.3cm}= G_\rmr (\ft,k)^{-1} - G_\rmr (\ft,k)^{-1} S_0^{\hlf} \Biggl[ - S_0   + \Biggl( - \frac{1}{2} I + \biggl( S_0^{\hlf} S_\ft S_0^{\hlf} \nonumber\\
   &\hspace{-0.1cm}+ \frac{1}{4} I  \biggr)^{\hlf} \Biggr)^{-1} + S_0^{\hlf} G_\rmr (\ft,k)^{-1} S_0^{\hlf}  \Biggr]^{-1} S_0^{\hlf} G_\rmr (\ft,k)^{-1} . \label{eq:Q_inv_backward}
\end{align}
Then by the same argument as for the time interval $ \bbra{k_\rmr, \ft} $, $ Q_{-,k} $ is also invertible on $ \bbra{0,k_\rmr -1} $.
Similarly, it can be shown that $ P_{-,k} $ is invertible on $ \bbra{0,\ft} $.

In the next step, we prove that $ I + B_k^\top Q_{-,k+1}^{-1} B_k \succ 0 $ for any $ k\in \bbra{0,\ft-1} $. Note that
\begin{align*}
	(I + B_k^\top Q_{-,k+1}^{-1} B_k)^{-1} &= I - B_k^\top (Q_{-,k+1} + B_k B_k^\top)^{-1} B_k \nonumber\\
   &= I - B_k^\top Q_{-,k}^{-1} B_k.
\end{align*}
Hence, it suffices to show $ I - B_k^\top Q_{-,k}^{-1} B_k\succ 0 $.
By \eqref{eq:Q_inv_forward}, for $ k \in \bbra{k_\rmr, \ft-1} $, $ \ft \ge 2 $, we have
\begin{align*}
	&I - B_k^\top Q_{-,k}^{-1} B_k = I + B_k^\top G_{\rmr} (k,0)^{-1} B_k \nonumber\\
	&\quad + B_k^\top G_{\rmr} (k,0)^{-1} S_0^{\hlf} \biggl[ S_0 + \frac{1}{2} I - \left( S_0^{\hlf} S_\ft S_0^{\hlf} + \frac{1}{4} I  \right)^{\hlf} \nonumber\\
   &\quad - S_0^{\hlf} G_\rmr (k,0)^{-1} S_0^{\hlf}  \biggr]^{-1} S_0^{\hlf} G_{\rmr} (k,0)^{-1} B_k .
\end{align*}
Since the expression in the square brackets is negative definite on $ \bbra{k_\rmr, \ft-1} $, it holds for sufficiently small $ \delta \in (0,1) $,
\begin{align}
	&I - B_k^\top Q_{-,k}^{-1} B_k \succ I + B_k^\top G_{\rmr} (k,0)^{-1} B_k  + B_k^\top G_{\rmr} (k,0)^{-1} S_0^{\hlf} \nonumber\\
   &\times  \left[  S_0 -\delta S_0 - S_0^{\hlf} G_\rmr (k,0)^{-1} S_0^{\hlf} \right]^{-1} S_0^{\hlf} G_\rmr (k,0)^{-1} B_k . \nonumber
\end{align}
Hence, we get
\begin{align}
	&I - B_k^\top Q_{-,k}^{-1} B_k \succ I + B_k^\top \left( G_\rmr (k,0) - (1-\delta)^{-1} I \right)^{-1} B_k \nonumber\\
	&= I - B_k^\top \left( \sum_{s=k}^{\ft-1} B_s B_s^\top + \left((1-\delta)^{-1} -1 \right) I \right)^{-1} B_k .
\end{align}
In addition, it holds that
\begin{align}
	&\Biggl[I - B_k^\top \left( \sum_{s=k}^{\ft-1} B_s B_s^\top + \left((1-\delta)^{-1} -1 \right) I \right)^{-1} B_k \Biggr]^{-1} \nonumber\\
	&= I + B_k^\top \left( \sum_{s=k+1}^{\ft-1} B_s B_s^\top + \left((1-\delta)^{-1} -1\right) I \right)^{-1} B_k \succ 0.\nonumber
\end{align}
Consequently, we obtain $ I + B_k^\top Q_{-,k+1}^{-1} B_k \succ 0$ for $ k\in \bbra{k_\rmr, \ft-1} $.
Noting also that $ Q_{-,k}^{-1} $ is given by \eqref{eq:Q_inv_backward} for $ k\in \bbra{0, k_\rmr-1} $, by the same argument above, it can be shown that $ I + B_k^\top Q_{-,k+1}^{-1} B_k \succ 0$ for $k\in \bbra{0,k_\rmr -1}$.

Next, we show the property (ii). Note that since $ 0\prec Q_{+,0} \prec I $, it holds $ Q_{+,\ft} \prec 0 $ by \eqref{eq:lyap_0N}. Thus by \eqref{eq:lyapunov_Qnew}, there exists $ s\in \bbra{0,\ft-1} $ such that $ Q_{+,s} \succ 0 $ and $ Q_{+,s+1} $ is not positive definite.
By assumption, $ Q_{+,s+1} $ is invertible and
\begin{align}
	&Q_{+,s+1}^{-1} =  (Q_{+,s} - B_s B_s^\top)^{-1}  \nonumber\\
   & = Q_{+,s}^{-1} + Q_{+,s}^{-1} B_s (I - B_s^\top Q_{+,s}^{-1} B_s)^{-1} B_s^\top Q_{+,s}^{-1} . \label{eq:Qs+1}
\end{align}
Now assume $ I - B_s^\top Q_{+,s}^{-1} B_s  \succ 0 $. Then by $ Q_{+,s}^{-1} \succ 0 $ and \eqref{eq:Qs+1}, we have $ Q_{+,s+1}^{-1} \succ 0 $, which contradicts the fact that $ Q_{+,s+1} $ is not positive definite.
Combining this with $ I + B_s^\top Q_{+,s+1}^{-1} B_s = (I - B_s^\top Q_{+,s}^{-1} B_s)^{-1} $, we conclude (ii).

Finally, by employing the relationship \eqref{eq:Q_Qnew_relation}, we obtain the desired result.

\section{Proofs of Proposition~\ref{prop:bridge_zero} and Theorem~\ref{thm:point_bridge}}\label{app:bridge_point}
\subsection{Proof of Proposition~\ref{prop:bridge_zero}}
	For simplicity, we write $ \sfP^\bullet_k := \sfP^\bullet (k,k), \ \wtilde{\sfP}_k :=\wtilde{\sfP} (k,k), \ \Phi_k := \Phi (\ft,k), \ G_{\rmr,k} := G_\rmr(\ft,k), \ G_{\rmc,k} := G_\rmc (\ft,k) $.
    First, we verify $ \sfP_k^\bullet = \wtilde{\sfP}_k $.
    Note that
	\begin{align*}
		&\sfP^\bullet_{k+1} = \what{A}_k \sfP^\bullet_k \what{A}_k^\top + B_k B_k^\top -  B_k B_k^\top \Phi_{k+1}^\top G_{\rmr,k}^\dagger \Phi_{k+1} B_k B_k^\top, \\
		&\sfP_0^\bullet = 0,
    \end{align*}
	and $ \wtilde{\sfP}_0 = 0 $. Thereafter, we show that $ \wtilde{\sfP}_k $ satisfies
	\begin{align}
		\wtilde{\sfP}_{k+1} &= \what{A}_k \wtilde{\sfP}_k \what{A}_k^\top + B_k B_k^\top -  B_k B_k^\top \Phi_{k+1}^\top G_{\rmr,k}^\dagger  \Phi_{k+1} B_k B_k^\top, \label{eq:Phat_transition}
	\end{align}
	which implies $ \sfP^\bullet_k = \wtilde{\sfP}_k $.

	Let $ T_k := \sfP_k + \sfQ_k $. Since $ T_{k+1} = A_k T_k A_k^\top $ by \eqref{eq:lyap_p_0} and \eqref{eq:lyap_q_0}, it holds $ T_k = \Phi (k,0) T_0 \Phi (k,0)^\top $. In addition, by the invertibility of $ A_k $ and $ T_\ft = \sfP_\ft = G_{\rmr,0} $, $ T_k $ is invertible for any $ k\in \bbra{0,\ft} $.
	Now, we have
	\begin{align*}
		T_k \Phi_{k}^\top \sfP_\ft^{-1} \Phi_{k} &= \Phi (k,0) T_0 \Phi (k,0)^\top \Phi_{k}^\top \sfP_\ft^{-1} \Phi_{k} \nonumber\\
		&= \Phi (k,\ft) \Phi (\ft,0) T_0 \Phi (\ft,0)^\top \sfP_\ft^{-1} \Phi_{k} \nonumber\\
		&= \Phi (k,\ft) T_\ft \sfP_\ft^{-1} \Phi (\ft,k) = I.
	\end{align*}
	Hence, it holds
	\begin{equation}\label{eq:PN}
		\Phi_{k}^\top \sfP_\ft^{-1} \Phi_{k} = (\sfP_k + \sfQ_k)^{-1} .
	\end{equation}
	By using this, we obtain
	\begin{align}
		&\wtilde{\sfP}_{k+1} =  \sfP_{k+1} - \sfP_{k+1} (\sfP_{k+1} + \sfQ_{k+1})^{-1} \sfP_{k+1} \nonumber\\
		&= A_k \sfP_k A_k^\top + B_k B_k^\top - (A_k \sfP_k A^\top + B_k B_k^\top) \nonumber\\
        &\quad \times  (A_k \sfP_k A_k^\top + A_k \sfQ_k A_k^\top)^{-1} (A_k \sfP_k A^\top + B_k B_k^\top) . \label{eq:Phat_next}
	\end{align}
	On the other hand, for the right-hand side of \eqref{eq:Phat_transition}, it holds
	\begin{align}
		&\what{A}_k \wtilde{\sfP}_k \what{A}_k^\top = (I - B_k B_k^\top \Phi_{k+1}^\top G_{\rmr,k}^\dagger \Phi_{k+1} ) \nonumber \\
		&\times \left(A_k \sfP_k A_k^\top - A_k \sfP_k A_k^\top (A_k \sfP_k A_k^\top + A_k \sfQ_k A_k^\top)^{-1} A_k \sfP_k A_k^\top \right) \nonumber\\
        &\times (I - B_k B_k^\top \Phi_{k+1}^\top G_{\rmr,k}^\dagger \Phi_{k+1} )^\top \nonumber\\
		&= A_k \sfP_k A_k^\top - A_k \sfP_k A_k^\top (A_k \sfP_k A_k^\top + A_k \sfQ_k A_k^\top)^{-1} A_k \sfP_k A_k^\top \nonumber \\
        &\quad - L_1 - L_1^\top + L_2, \label{eq:AhatPhatAhat}
	\end{align}
    where
    \begin{align*}
        &L_1 :=  B_k B_k^\top \Phi_{k+1}^\top G_{\rmr,k}^\dagger \Phi_{k+1} A_k \sfP_k A_k^\top \nonumber\\
        &\qquad\quad \times \left(I - (A_k \sfP_k A_k^\top + A_k \sfQ_k A_k^\top)^{-1} A_k \sfP_k A_k^\top \right), \\
		&L_2 := B_k B_k^\top \Phi_{k+1}^\top G_{\rmr,k}^\dagger \Phi_{k+1} \nonumber \\ 
        &\times \left(A_k \sfP_k A_k^\top - A_k \sfP_k A_k^\top (A_k \sfP_k A_k^\top + A_k \sfQ_k A_k^\top)^{-1} A_k \sfP_k A_k^\top \right)  \nonumber\\
		&\times \Phi_{k+1}^\top G_{\rmr,k}^\dagger \Phi_{k+1} B_k B_k^\top .
    \end{align*}
	Here, we have
	\begin{align*}
		&L_1 =  B_k B_k^\top \Phi_{k+1}^\top G_{\rmr,k}^\dagger \Phi_{k+1}  ( A_k \sfP_k A_k^\top + A_k \sfQ_k A_k^\top) \nonumber\\
        &\qquad \times \left(I - ( A_k \sfP_k A_k^\top + A_k\sfQ_k A_k^\top)^{-1} A_k \sfP_k A_k^\top\right) \nonumber\\
		&\qquad - B_k B_k^\top \Phi_{k+1}^\top G_{\rmr,k}^\dagger \Phi_{k+1} A_k \sfQ_k A_k^\top \nonumber\\
        &\qquad \times \left(I - ( A_k\sfP_k A_k^\top + A_k\sfQ_k A_k^\top)^{-1} A_k\sfP_k A_k^\top \right) \nonumber\\
		&= B_kB_k^\top \Phi_{k+1}^\top G_{\rmr,k}^\dagger \Phi_{k+1} A_k \sfQ_k A_k^\top ( A_k\sfP_k A_k^\top + A_k\sfQ_k A_k^\top)^{-1} \nonumber\\
		&\quad \times A_k\sfP_k A_k^\top .
	\end{align*}
	Moreover, it follows from
	\begin{equation}\label{eq:Q_gram}
		A_k \sfQ_k A_k^\top = \Phi(k+1,\ft) G_{\rmr,k} \Phi(k+1,\ft)^\top
	\end{equation}
	that
	\begin{align*}
		L_1 &= B_k B_k^\top \Phi_{k+1}^\top G_{\rmr,k}^\dagger G_{\rmr,k} \Phi_{k+1}^{-\top} \nonumber\\
        &\quad \times (A_k \sfP_k A_k^\top + A_k \sfQ_k A_k^\top)^{-1} A_k \sfP_k A_k^\top .
	\end{align*}
	It is known that if $ {\rm ran} (Y^\top) \supseteq {\rm ran} (X^\top) $, then $ Y = Y X^\dagger X $~\cite[Theorem~3.11]{Yanai2011}, and we have
	\begin{align*}
		&{\rm ran} (G_{\rmr,k}) = {\rm ran} ( [\Phi_{k+1} B_k ~~ \Phi_{k+2} B_{k+1}  \ \cdots \ B_{\ft-1}] ) \nonumber\\
        & \supseteq {\rm ran} (\Phi_{k+1} B_k) = {\rm ran} ( \Phi_{k+1} B_k B_k^\top \Phi_{k+1}^\top ) .
	\end{align*}
    Hence, we have
    \begin{equation}\label{eq:cancel}
        B_k B_k^\top \Phi_{k+1}^\top G_{\rmr,k}^\dagger G_{\rmr,k} = B_k B_k^\top \Phi_{k+1}^\top,
		\end{equation}
	which yields
	\begin{align}\label{eq:J1}
		L_1 = B_k B_k^\top (A_k \sfP_k A_k^\top + A_k \sfQ_k A_k^\top)^{-1} A_k \sfP_k A_k^\top .
	\end{align}
	By a similar argument, 
	\begin{align}
		&L_2 = B_kB_k^\top \Phi_{k+1}^\top G_{\rmr,k}^\dagger \Phi_{k+1} A_k\sfQ_{k} A_k^\top \nonumber\\
		& \times (A_k\sfP_k A_k^\top + A_k\sfQ_k A_k^\top)^{-1} A_k\sfP_k A_k^\top \Phi_{k+1}^\top G_{\rmr,k}^\dagger \Phi_{k+1} B_kB_k^\top \nonumber\\
		&=B_kB_k^\top (A_k\sfP_k A_k^\top + A_k\sfQ_k A_k^\top)^{-1} \nonumber\\
		&\quad \times A_k\sfP_k A_k^\top \Phi_{k+1}^\top G_{\rmr,k}^\dagger \Phi_{k+1} B_kB_k^\top \nonumber\\
		&=B_kB_k^\top (A_k\sfP_k A_k^\top + A_k\sfQ_k A_k^\top)^{-1} (A_k\sfP_k A_k^\top + A_k\sfQ_k A_k^\top) \nonumber\\
        &\quad \times \Phi_{k+1}^\top G_{\rmr,k}^\dagger \Phi_{k+1} B_kB_k^\top \nonumber\\
        &\quad  - B_kB_k^\top (A_k\sfP_k A_k^\top + A_k\sfQ_k A_k^\top)^{-1}  A_k\sfQ_k A_k^\top  \nonumber\\
				&\quad \times \Phi_{k+1}^\top G_{\rmr,k}^\dagger \Phi_{k+1} B_kB_k^\top \nonumber\\
		&= B_kB_k^\top \Phi_{k+1}^\top G_{\rmr,k}^\dagger \Phi_{k+1} B_kB_k^\top \nonumber\\
        &\quad - B_kB_k^\top (A_k\sfP_k A_k^\top + A_k\sfQ_k A_k^\top)^{-1} B_kB_k^\top . \label{eq:J2}
	\end{align}
	By \eqref{eq:Phat_next},~\eqref{eq:AhatPhatAhat},~\eqref{eq:J1},~\eqref{eq:J2}, we obtain \eqref{eq:Phat_transition}, and thus $ \sfP^\bullet_k = \wtilde{\sfP}_k $.

    Next, we show $ \sfP^\bullet (k,s) = \wtilde{\sfP}(k,s) $ for $ k\le s $. Denote by $ \what{\Phi} (k,l) $ the state-transition matrix for $ \what{A}_k $ as in \eqref{eq:transition_matrix}.
    Then,
	\begin{align*}
		\sfP^\bullet (k,s) &= \bbE\left[x_k^\bullet \biggl(\what{\Phi} (s,k) x_k^\bullet + \sum_{l=k}^{s-1} \what{\Phi} (s,l+1) B_l w_l^\bullet\biggr)^\top \right] \nonumber \\
        &= \sfP^\bullet (k,k) \what{\Phi} (s,k)^\top, \ k\le s .
	\end{align*}
	Therefore, it holds
	\begin{equation}\label{eq:Pbullet_stransition}
		\sfP^\bullet (k,s+1) = \sfP^\bullet (k,s) \what{A}_s^\top, \ k\le s.
	\end{equation}
    In what follows, we show that the same relationship
    \begin{equation}\label{eq:Phat_stransition}
        \wtilde{\sfP}(k,s+1) = \wtilde{\sfP}(k,s) \what{A}_s^\top, \ k\le s
    \end{equation}
    also holds. The left-hand side of \eqref{eq:Phat_stransition} is written as
	\begin{align}
		&\wtilde{\sfP}(k,s+1) = \sfP_k \Phi (s+1,k)^\top- \sfP_k \Phi_{k}^\top \sfP_\ft^{-1} \Phi_{s} \sfP_s A_s^\top \nonumber\\
        &\qquad\qquad\qquad - \sfP_k \Phi_{k}^\top \sfP_\ft^{-1} \Phi_{s+1} B_s B_s^\top \nonumber\\
		&=\sfP_k \Phi (s+1,k)^\top - \sfP_k \Phi_{k}^\top \sfP_\ft^{-1} \Phi_{s} \sfP_s A_s^\top \nonumber\\
        &\quad - \sfP_k \Phi (s+1,k)^\top (A_s\sfP_s A_s^\top + A_s\sfQ_s A_s^\top)^{-1} B_s B_s^\top , \label{eq:Phat_snext}
	\end{align}
	where in the last line we used \eqref{eq:PN}.
	On the other hand, the right-hand side of \eqref{eq:Phat_stransition} is
	\begin{align}
		&\wtilde{\sfP}(k,s) \what{A}_s^\top  = \left(\sfP_k \Phi(s,k)^\top - \sfP_k \Phi_{k}^\top \sfP_\ft^{-1} \Phi_{s} \sfP_{s} \right) A_s^\top \nonumber\\
		&\qquad\qquad\qquad \times (I - \Phi_{s+1}^\top G_{\rmr,s}^\dagger \Phi_{s+1} B_s B_s^\top) \nonumber\\
		&= \sfP_k \Phi (s+1,k)^\top - \sfP_k \Phi_{k}^\top \sfP_\ft^{-1} \Phi_{s} \sfP_s A_s^\top \nonumber\\
		&\quad - \sfP_k \Phi_{k}^\top G_{\rmr,s}^\dagger \Phi_{s+1} B_s B_s^\top + L_3, \label{eq:Phat_Ahat}
	\end{align}
    where
    \begin{align*}
        L_3 &:= \sfP_k \Phi(s+1,k)^\top (A_s\sfP_sA_s^\top + A_s\sfQ_s A_s^\top)^{-1} A_s\sfP_s A_s^\top \\
        &\quad \times \Phi_{s+1}^\top G_{\rmr,s}^\dagger \Phi_{s+1} B_s B_s^\top \\
        &= \sfP_k \Phi (s+1,k)^\top \Phi_{s+1}^\top G_{\rmr,s}^\dagger \Phi_{s+1} B_s B_s^\top \\
		&\quad - \sfP_k \Phi (s+1,k)^\top (A_s \sfP_sA_s^\top + A_s\sfQ_s A_s^\top)^{-1} A_s\sfQ_s A_s^\top \\
        &\quad \times \Phi_{s+1}^\top G_{\rmr,s}^\dagger \Phi_{s+1} B_s B_s^\top \\
		&= \sfP_k \Phi_{k}^\top G_{\rmr,s}^\dagger \Phi_{s+1} B_s B_s^\top \\
        &\quad  - \sfP_k \Phi (s+1,k)^\top (A_s \sfP_s A_s^\top + A_s\sfQ_s A_s^\top)^{-1} B_s B_s^\top .
    \end{align*}
    By \eqref{eq:Phat_snext} and \eqref{eq:Phat_Ahat}, we obtain \eqref{eq:Phat_stransition}.
    Lastly, by combining \eqref{eq:Pbullet_stransition},~\eqref{eq:Phat_stransition} with $ \sfP_k^\bullet = \wtilde{\sfP}_k $, we arrive at $ \sfP^\bullet (k,s) = \wtilde{\sfP} (k,s) $ for $ k\le s $, which completes the proof.

	\subsection{Proof of Theorem~\ref{thm:point_bridge}}
  Here we show \eqref{eq:first_moment}.
	First, by using \eqref{eq:PN}, we obtain
	\begin{align}
		\ell_{k+1} &=\Phi (k+1,0) \bar{x}_0 + (A_k \sfP_k A_k^\top + B_k B_k^\top) A_k^{-\top} \nonumber\\
        &\quad \times (\sfP_k + \sfQ_k)^{-1} \Phi (k,\ft) (\bar{x}_\ft - \Phi(\ft,0) \bar{x}_0) \nonumber\\
		&= \Phi (k+1,0) \bar{x}_0  \nonumber\\
        &\quad +A_k \sfP_k (\sfP_k + \sfQ_k)^{-1} \Phi (k,\ft) (\bar{x}_\ft - \Phi(\ft,0) \bar{x}_0) \nonumber\\
		&\quad + B_k B_k^\top (A_k \sfP_k A_k^\top + A_k \sfQ_k A_k^\top)^{-1} \Phi (k+1,\ft) \nonumber\\ 
        &\quad \times (\bar{x}_\ft - \Phi(\ft,0) \bar{x}_0) . \label{eq:l_next}
	\end{align}
	On the other hand, for the right-hand side of \eqref{eq:first_moment}, it holds
	\begin{align}
		\what{A}_k \ell_k &= \Phi (k+1,0) \bar{x}_0 - B_k B_k^\top \Phi (\ft,k+1)^\top G_\rmr (\ft,k)^\dagger \nonumber \\
        &\quad \times \Phi (\ft,0) \bar{x}_0 + L_4 + L_5, \label{eq:Ahat_l}
	\end{align}
    where
    \begin{align*}
        L_4 &:= A_k \sfP_k \Phi (\ft,k)^\top \sfP_\ft^{-1} (\bar{x}_\ft - \Phi (\ft,0) \bar{x}_0), \\
		L_5 &:= - B_k B_k^\top \Phi (\ft,k+1)^\top G_\rmr (\ft,k)^\dagger \Phi (\ft,k+1) \nonumber\\
        &\quad \times A_k \sfP_k (\sfP_k + \sfQ_k)^{-1} \Phi (k,\ft) (\bar{x}_\ft - \Phi(\ft,0) \bar{x}_0) .
    \end{align*}
	Here, we have
	\begin{align}
		&L_4 = A_k \sfP_k (\sfP_k + \sfQ_k)^{-1} \Phi (k,\ft) (\bar{x}_\ft - \Phi(\ft,0) \bar{x}_0), \label{eq:J4}\\
		&L_5 = - B_k B_k^\top \Phi (\ft,k+1)^\top G_\rmr (\ft,k)^\dagger \Phi (\ft,k+1) A_k \sfP_k A_k^\top \nonumber\\
        &\times  (A_k \sfP_k A_k^\top  +  A_k \sfQ_k A_k^\top)^{-1} \Phi (k+1,\ft) (\bar{x}_\ft - \Phi(\ft,0) \bar{x}_0)\nonumber\\
		&= - B_k B_k^\top \Phi(\ft,k+1)^\top G_\rmr (\ft,k)^\dagger (\bar{x}_\ft - \Phi(\ft,0) \bar{x}_0) \nonumber\\
		& + B_k B_k^\top \Phi (\ft,k+1)^\top G_\rmr (\ft,k)^\dagger \Phi(\ft,k+1) A_k \sfQ_k A_k^\top \nonumber\\
        &\times (A_k \sfP_k A_k^\top + A_k \sfQ_k A_k^\top)^{-1} \Phi(k+1,\ft) (\bar{x}_\ft - \Phi(\ft,0) \bar{x}_0) \nonumber\\
		&=- B_k B_k^\top \Phi(\ft,k+1)^\top G_\rmr (\ft,k)^\dagger (\bar{x}_\ft - \Phi(\ft,0) \bar{x}_0) \nonumber\\
		&\quad + B_k B_k^\top (A_k \sfP_k A_k^\top + A_k \sfQ_k A_k^\top)^{-1} \Phi(k+1,\ft) \nonumber\\
        &\quad \times(\bar{x}_\ft - \Phi(\ft,0) \bar{x}_0) , \label{eq:J5}
	\end{align}
	where in the last line we used \eqref{eq:Q_gram},~\eqref{eq:cancel}.
  By \eqref{eq:l_next}, \eqref{eq:Ahat_l}, \eqref{eq:J4}, \eqref{eq:J5}, we get \eqref{eq:first_moment}.

\section{Proof of Theorem~\ref{thm:schrodinger}}\label{app:schrodinger_proof}
By the expression \eqref{eq:separate_KL}, it suffices to prove the following two statements:
	\begin{itemize}
		\item[(i)] Consider the problem of minimizing $ D_{\rm KL}(\what{\bbP}_{0,\ft} \| \wtilde{\bbP}_{0,\ft} ) $ with respect to $ \what{\bbP}_{0,\ft} $ under the constraint that the marginals of $ \what{\bbP}_{0,\ft} $ at $ k=0 $ and $ k = \ft $ are equal to $ \calN(0,\bar{\Sigma}_0) $ and $\calN(0,\bar{\Sigma}_\ft) $, respectively. Then, the joint distribution $ \bbP_{0,\ft}^* $ of $ (x_0^*, x_\ft^*) $ is an optimal solution to this problem;
		\item[(ii)] The conditional distributions of $ \{x_k^*\} $ and $ \{\wtilde{x}_k\} $ satisfy
		\begin{align}
			&\bbP^*_{1:\ft-1|0,\ft} (\cdot | x_0^* = \bar{x}_0, x_\ft^* = \bar{x}_\ft) \nonumber\\
			&= \wtilde{\bbP}_{1:\ft-1|0,\ft} (\cdot | \wtilde{x}_0 = \bar{x}_0, \wtilde{x}_\ft = \bar{x}_\ft), \  \forall \bar{x}_0, \bar{x}_\ft \in \bbR^n .
		\end{align} 
	\end{itemize}

	{\bf Proof of part~(i).} For brevity, we denote $ \calN_\Sigma (\cdot) := \calN (\cdot | 0, \Sigma) $.
    Note that for any $ \Sigma, \Xi \succ 0 $, the KL divergence between $ n $-variate Gaussian distributions $ \calN_\Sigma , \calN_{\Xi} $ is given by
	\begin{align}
		&\kl{\calN_\Sigma}{\calN_{\Xi}} \nonumber\\
		&= \frac{1}{2} \left( \log \det (\Xi) -  \log\det (\Sigma)  +  \trace (\Xi^{-1} \Sigma) -  n \right) .\label{eq:kl_gauss}
	\end{align}
    Now, let $ p_\Sigma $ be the density function of a distribution with covariance matrix $ \Sigma $, which is not necessarily Gaussian. 
		Then, it is known that given covariance matrices $ \Sigma $ and $ \Xi $, the density $ p_\Sigma $ that minimizes $ \kl{p_\Sigma}{\calN_\Xi} $ is the Gaussian distribution $ \calN_\Sigma $; see the proof of \cite[Theorem~11]{Chen2016part1}.

	Now, we introduce
	\begin{align*}
		&\Sigma_{0,\ft} := \begin{bmatrix}
			\bar{\Sigma}_0 & Y^\top\\ Y & \bar{\Sigma}_\ft
		\end{bmatrix}, \ Y \in \bbR^{n\times n} \\
		&\Xi_{0,\ft} := \begin{bmatrix}
			\Xi_0 & \bar{\Sigma}_0 \Phi (\ft,0)^\top \\
			\Phi (\ft,0) \bar{\Sigma}_0 & \Xi_\ft
		\end{bmatrix}, \\
		&\Xi_0 := \bsigma_0, \ \ \Xi_k := \Phi(k,0) \bsigma_0 \Phi(k,0)^\top + G_\rmr (k,0), \ k\in \bbra{1,\ft}.
	\end{align*}
	Then $ \Xi_{0,\ft} $ is the covariance of the Gaussian distribution $ \wtilde{\bbP}_{0,\ft} $.
    In the light of the above description, the minimization of $ D_{\rm KL} (\what{\bbP}_{0,\ft}\|\wtilde{\bbP}_{0,\ft}) $ is equivalent to the minimization of $ \kl{\calN_{\Sigma_{0,\ft}}}{\calN_{\Xi_{0,\ft}}} $ with respect to $ Y $.
    By the formula \eqref{eq:kl_gauss}, this is further equivalent to the minimization of
	\begin{equation}\label{eq:min_sub}
		\trace ( \Xi_{0,\ft}^{-1} \Sigma_{0,\ft} ) - \log \det (\Sigma_{0,\ft})
	\end{equation}
    with respect to $ Y $.
	Noting that
	\begin{equation*}
		\Xi_{0,\ft} = \begin{bmatrix}
			I \\ \Phi (\ft,0)
		\end{bmatrix} \bar{\Sigma}_0
		\begin{bmatrix}
			I & \Phi (\ft,0)^\top
		\end{bmatrix} +
		\begin{bmatrix}
			0 & 0\\ 0 & G_\rmr (\ft,0)
		\end{bmatrix},
	\end{equation*}
	we have
	\begin{equation*}
		\Xi_{0,\ft}^{-1} =
		\begin{bmatrix}
			\bar{\Sigma}_0^{-1} + \Phi_0^\top G_{\rmr,0}^{-1} \Phi_0 & - \Phi_0^\top G_{\rmr,0}^{-1} \\
			- G_{\rmr,0}^{-1} \Phi_0 & G_{\rmr,0}^{-1} 
		\end{bmatrix} ,
	\end{equation*}
	where we introduced the abbreviations $ \Phi_0 := \Phi (\ft,0) $, $G_{\rmr,0} := G_\rmr (\ft,0) $.
    Then, \eqref{eq:min_sub} is written as
	\begin{align*}
		&\trace\Bigl( (\bar{\Sigma}_0^{-1} + \Phi_0^\top G_{\rmr,0}^{-1} \Phi_0) \bar{\Sigma}_0 - \Phi_0^\top G_{\rmr,0}^{-1} Y - Y^\top G_{\rmr,0}^{-1} \Phi_0 \\
		&+ G_{\rmr,0}^{-1} \bar{\Sigma}_\ft  \Bigr) - \log \det (\bar{\Sigma}_0) - \log \det (\bar{\Sigma}_\ft - Y \bar{\Sigma}_0^{-1} Y^\top) .
	\end{align*}

	In what follows, we consider maximizing
	\begin{equation}
		f(Y) := \log \det (\bar{\Sigma}_\ft - Y \bar{\Sigma}_0^{-1} Y^\top) + 2 \trace (\Phi^\top G_{\rmr,0}^{-1} Y),
	\end{equation}
    which is concave. By taking the derivative of $ f $, the necessary and sufficient condition for $ Y $ to be the maximizer of $ f(Y) $ is obtained as
	\begin{equation}\label{eq:first_order}
		- 2 \bar{\Sigma}_0^{-1} Y^\top (\bar{\Sigma}_\ft - Y \bar{\Sigma}_0^{-1} Y^\top )^{-1} + 2\Phi_0^\top G_{\rmr,0}^{-1} = 0 .
	\end{equation}
	For notational simplicity, we omit the subscript ``$ - $'' of $ Q_{-,k} $.
	Define the state-transition matrix $ \Phi_{\qm} (k,l) $ for $ A_{Q,k} $ like in \eqref{eq:transition_matrix}.
	Then $ \bbP_{0,\ft}^* $ is Gaussian with covariance matrix
	\begin{equation*}
		\begin{bmatrix}
			\bsigma_0 & \bar{\Sigma}_0 \Phi_\qm (\ft,0)^\top \\
			\Phi_\qm (\ft,0) \bar{\Sigma}_0 & \bsigma_\ft
		\end{bmatrix} .
	\end{equation*}
    Hence, it suffices to show that $ Y = \Phi_\qm (\ft,0) \bsigma_0 $ satisfies the condition \eqref{eq:first_order}, which is equivalent to
	\begin{align}
		&\Phi_\qm (\ft,0)^\top \left(\bsigma_\ft - \Phi_\qm (\ft,0) \bsigma_0 \Phi_\qm (\ft,0)^\top \right)^{-1} \nonumber\\
        &= \Phi (\ft,0)^\top G_\rmr (\ft,0)^{-1} \nonumber\\
		&= \Phi (\ft,0)^\top (\Xi_\ft - \Phi (\ft,0) \bar{\Sigma}_0 \Phi(\ft,0)^\top)^{-1} . \label{eq:Y_optimality}
	\end{align}
    Taking the inverse of both sides of \eqref{eq:Y_optimality} yields
	\begin{equation*}
		\bsigma_\ft \Phi_\qm (0,\ft)^\top - \Phi_\qm (\ft,0) \bsigma_0 = \Xi_\ft \Phi (0,\ft)^\top - \Phi (\ft,0) \bsigma_0 .
	\end{equation*}
	Hereafter, we show that
	\begin{align}\label{eq:max_condition}
		&\Sigma_k \Phi_\qm (0,k)^\top - \Phi_\qm (k,0) \bsigma_0 = \Xi_k \Phi (0,k)^\top - \Phi (k,0) \bsigma_0, \nonumber\\
        &\hspace{5.5cm}\forall k\in \bbra{0,\ft} ,
	\end{align}
    where $ \Sigma_k = \bbE[x_k^* (x_k^*)^\top] $, that is,
	\[
		\Sigma_{k+1} = A_{Q,k} \Sigma_k A_{Q,k}^\top + B_k (I + B_k^\top \qmk{k+1}^{-1} B_k)^{-1} B_k^\top.
	\]
    For this purpose, we define
	\begin{align*}
		&F_{1,k} := \Sigma_k \Phi_\qm (0,k)^\top - \Phi_\qm (k,0) \bsigma_0, \\
		&F_{2,k} := \Xi_k \Phi (0,k)^\top - \Phi (k,0) \bsigma_0, \\
		&F_{3,k} := \qmk{k} (\Phi_\qm (0,k)^\top - \Phi (0,k)^\top),
	\end{align*}
    and show $ F_{2,k} = F_{3,k} $ and $ F_{1,k} = F_{3,k} $, which imply \eqref{eq:max_condition}.

    First, noting that $ F_{2,0} = F_{3,0} $, we check that $ F_{2,k} $ and $ F_{3,k} $ satisfy the same difference equation, which means $ F_{2,k} = F_{3,k} $ for all $ k\in \bbra{0,\ft} $. For $ F_{2,k} $, we have
	\begin{align*}
		F_{2,k+1} &= \left( \Phi(k+1,0) \bsigma_0 \Phi(k+1,0)^\top + G_\rmr (k+1,0) \right) \nonumber\\
        &\quad\times \Phi (0,k+1)^\top - \Phi (k+1,0) \bsigma_0 \nonumber\\
		&= G_\rmr (k+1,0) \Phi (0,k+1)^\top \nonumber\\
		&= \left( A_k G_\rmr (k,0) A_k^\top + B_kB_k^\top  \right) \Phi (0,k+1)^\top \nonumber\\
		&= A_k F_{2,k} + B_k B_k^\top \Phi(0,k+1)^\top .
	\end{align*}
	On the other hand, for $ F_{3,k} $, it holds
	\begin{align}
		&F_{3,k+1} = \qmk{k+1} \Phi_\qm (0,k+1)^\top - (A_k \qmk{k} A_k^\top - B_k B_k^\top) \nonumber\\
        &\qquad\qquad \times\Phi(0,k+1)^\top  \nonumber\\
		&= (\qmk{k+1} + B_k B_k^\top) A_k^{-\top} \Phi_{\qm} (0,k)^\top - A_k \qmk{k} \Phi (0,k)^\top \nonumber\\
        &\quad + B_k B_k^\top \Phi (0,k+1)^\top  \nonumber\\
		&= A_k F_{3,k} + B_k B_k^\top \Phi (0,k+1)^\top ,
	\end{align}
	where we used $ A_{Q,k} = (I + B_k B_k^\top \qmk{k+1}^{-1})^{-1} A_k $.
    Hence, we obtain $ F_{2,k} = F_{3,k} $.

    Next, define $ M_k := \qmk{k}^{-1} (F_{3,k} - F_{1,k}) $. Then
	\begin{align}
		M_{k+1} &= (I - \qmk{k+1}^{-1} \Sigma_{k+1}) \Phi_\qm (0,k+1)^\top \nonumber\\
        &\quad + \qmk{k+1}^{-1} \Phi_\qm (k+1,0) \bsigma_0 - \Phi (0,k+1)^\top . \label{eq:M_transition}
	\end{align}
	The first term is written as
	\begin{align}
		&(I - \qmk{k+1}^{-1} \Sigma_{k+1})\Phi_\qm (0,k+1)^\top \nonumber\\
		&= \Bigl( A_{Q,k}^{-\top} -  \bigl[ (\qmk{k+1} + B_k B_k^\top )^{-1} A_k \Sigma_k \nonumber\\
        & \quad + \qmk{k+1}^{-1} B_k (I + B_k^\top \qmk{k+1}^{-1} B_k)^{-1} B_k^\top A_{Q,k}^{-\top}  \bigr]  \Bigr) \Phi_\qm (0,k)^\top . \label{eq:M_first_term}
	\end{align}
	Here, noting that
    \begin{align*}
        A_{Q,k}^{-\top} &= \left( \qmk{k+1} - B_k (I + B_k^\top \qmk{k+1}^{-1} B_k)^{-1} B_k^\top  \right)^{-1} \qmk{k+1} A_k^{-\top},
    \end{align*}
    we have
	\begin{align*}
		&\qmk{k+1}^{-1} B_k(I + B_k^\top \qmk{k+1}^{-1} B_k)^{-1} B_k^\top A_{Q,k}^{-1} 
		= - A_k^{-\top} + A_{Q,k}^{-\top} .
	\end{align*}
	Thus, \eqref{eq:M_first_term} is written as
	\begin{align}
		&(I - \qmk{k+1}^{-1} \Sigma_{k+1})\Phi_\qm (0,k+1)^\top \nonumber\\
		&= \left( - (A_k \qmk{k} A_k^\top)^{-1} A_k \Sigma_k + A_k^{-\top}     \right) \Phi_\qm (0,k)^\top \nonumber\\
		&= A_k^{-\top} ( I -\qmk{k}^{-1} \Sigma_k ) \Phi_\qm (0,k)^\top . \label{eq:M_first_term2}
	\end{align}
	For the second term of the right-hand side of \eqref{eq:M_transition}, we have
	\begin{align}
		\qmk{k+1}^{-1} \Phi_\qm (k+1,0) \bar{\Sigma}_0
		&= (\qmk{k+1} + B_k B_k^\top)^{-1} A_k \Phi_\qm (k,0) \bar{\Sigma}_0 \nonumber\\
		&= A_k^{-\top} \qmk{k}^{-1} \Phi_\qm (k,0) \bsigma_0 . \label{eq:M_second_term}
	\end{align}
	It follows from \eqref{eq:M_transition},~\eqref{eq:M_first_term2},~\eqref{eq:M_second_term} that
	\[
		M_{k+1} = A_k^{-\top} M_k .
	\]
    Combining this with $ M_0 = 0 $, we obtain $ M_k = 0 $ for all $ k\in \bbra{0,\ft} $, that is, $ F_{1,k} = F_{3,k}$ for all $ k\in \bbra{0,\ft} $.
	In summary, we get \eqref{eq:max_condition}, which implies (i).

	{\bf Proof of part~(ii).} By Theorem~\ref{thm:point_bridge}, it suffices to show that for any end points $ \bar{x}_0, \bar{x}_\ft \in \bbR^n$, the pinned process associated to \eqref{eq:noise_driven_point} coincides with the one associated to \eqref{eq:opt_transition}.
    In other words, we verify that for any $ k\in \bbra{0,\ft-1} $,
	\begin{align}
		&B_k B_k^\top \Phi (\ft,k+1)^\top G_\rmr (\ft,k)^\dagger  \nonumber\\
        &\quad = B_{Q,k} B_{Q,k}^\top \Phi_{\qm} (\ft,k+1)^\top G_{Q,\rmr} (\ft,k)^\dagger, \label{eq:final_coef}\\
		&\what{A}_k = \bigl( I - B_{Q,k} B_{Q,k}^\top \Phi_\qm (\ft,k+1)^\top G_{Q,\rmr} (\ft,k)^\dagger \nonumber\\
        &\qquad \times \Phi_\qm (\ft,k+1)  \bigr) A_{Q,k} =: \what{A}_{Q,k} , \label{eq:bridge_A} \\
		&\Lambda_k := B_k B_k^\top - B_k B_k^\top \Phi (\ft,k+1)^\top G_\rmr (\ft,k)^\dagger \nonumber\\
        &\qquad \times \Phi (\ft,k+1)B_k B_k^\top \nonumber\\
		&\quad = B_{Q,k} B_{Q,k}^\top - B_{Q,k} B_{Q,k}^\top \Phi_\qm (\ft,k+1)^\top G_{Q,\rmr} (\ft,k)^\dagger \nonumber\\
        &\qquad \times \Phi_\qm (\ft,k+1) B_{Q,k} B_{Q,k}^\top =: \Lambda_{Q,k} . \label{eq:bridge_cov}
	\end{align}
    Here, $ G_{Q,\rmr} $ and $ G_{Q,\rmc} $ are the reachability Gramian and the controllability Gramian where $ \Phi $ and $ B_k $ in \eqref{eq:reachability_gramian} and \eqref{eq:controllability_gramian} are replaced by $ \Phi_{\qm} $ and $ B_{Q,k} $, respectively.

	First, note that
	\begin{align}
		&\Phi (\ft,k+1) B_k B_k^\top \Phi (\ft,k+1)^\top G_\rmr (\ft,k)^\dagger \nonumber\\
        &\quad= \Phi (\ft,k+1) B_k B_k^\top A_k^{-\top} \left( G_\rmr (\ft,k) \Phi (k,\ft)^\top \right)^\dagger .\label{eq:Gr_dagger}
	\end{align}
    Indeed, since it holds $ {\rm ran} (G_\rmr (\ft,k)) \supseteq {\rm ran} (\Phi (\ft,k+1) B_k B_k^\top \Phi (\ft,k+1)^\top ) $, there exists $ W \in \bbR^{n\times n} $ such that
	\begin{align*}
		\Phi (\ft,k+1) B_k B_k^\top \Phi (\ft,k+1)^\top &= G_\rmr (\ft,k) W \\
		&= W^\top G_\rmr (\ft,k),
	\end{align*}
	and we have
	\begin{align*}
		&\Phi (\ft,k+1) B_k B_k^\top \Phi (\ft,k+1)^\top G_\rmr (\ft,k)^\dagger \nonumber\\
        &\quad = W^\top G_\rmr (\ft,k) G_\rmr (\ft,k)^\dagger, \\
		&\Phi (\ft,k+1) B_k B_k^\top A_k^{-\top} \left( G_\rmr (\ft,k) \Phi (k,\ft)^\top \right)^\dagger \nonumber\\
        &\quad = W^\top G_\rmr (\ft,k) \Phi (k,\ft)^\top (G_\rmr (\ft,k) \Phi (k,\ft)^\top)^\dagger .
	\end{align*}
    Here, $ G_\rmr (\ft,k) G_\rmr (\ft,k)^\dagger $ is the orthogonal projection matrix onto $ {\rm ran} (G_\rmr (\ft,k)) $, and $ G_\rmr (\ft,k) \Phi (k,\ft)^\top (G_\rmr (\ft,k) \Phi (k,\ft)^\top)^\dagger $ is also the orthogonal projection onto $ {\rm ran} (G_\rmr (\ft,k) \Phi (k,\ft)^\top) = {\rm ran} (G_\rmr (\ft,k) ) $.
    Thus, by the uniqueness of the orthogonal projection, we obtain \eqref{eq:Gr_dagger}.

    By using \eqref{eq:control_reachable},~\eqref{eq:Gr_dagger}, the left-hand side of \eqref{eq:final_coef} is equal to
	\begin{align}
         B_k B_k^\top A_k^{-\top} \left( \Phi (\ft,k) G_\rmc (\ft,k) \right)^\dagger. \label{eq:left_final_coeff}
	\end{align}
	Similarly, the right-hand side of \eqref{eq:final_coef} is written as
	\begin{align}
		&B_{Q,k} B_{Q,k}^\top A_{Q,k}^{-\top} \left( \Phi_\qm (\ft,k) G_{Q,\rmc} (\ft,k)  \right)^\dagger \nonumber\\
        & = B_k B_k^\top A_k^{-\top} \left( \Phi_\qm (\ft,k) G_{Q,\rmc} (\ft,k)  \right)^\dagger , \label{eq:right_final_coeff}
	\end{align}
	where we used $ B_{Q,k} B_{Q,k}^\top A_{Q,k}^{-\top} = B_k B_k^\top A_k^{-\top} $.
    For notational simplicity, we write $ R_{1,k} := G_\rmc (\ft,k), R_{2,k} := G_{Q,\rmc} (\ft,k) $ and set $ R_{1,\ft} = R_{2,\ft} = 0 $. By \eqref{eq:left_final_coeff} and \eqref{eq:right_final_coeff}, to verify \eqref{eq:final_coef}, it suffices to check that
	\begin{equation}\label{eq:R1_R2}
		\Phi (\ft,k) R_{1,k} - \Phi_\qm (\ft,k) R_{2,k} = 0 .
	\end{equation}
	Since it holds
	\begin{align*}
		A_{Q,k} = (I + B_k B_k^\top \qmk{k+1}^{-1})^{-1} A_k = \qmk{k+1} A_k^{-\top} \qmk{k}^{-1},
	\end{align*}
	we have $ \Phi_\qm (\ft,k)  = \qmk{\ft} \Phi (k,\ft)^\top \qmk{k}^{-1} $, and therefore 
	\begin{align*}
		&\Phi (\ft,k) R_{1,k} - \Phi_\qm (\ft,k) R_{2,k} \nonumber\\
        &= \Phi (\ft,k) \left( R_{1,k} - \Phi (k,\ft) \qmk{\ft} \Phi (k,\ft)^\top \qmk{k}^{-1} R_{2,k}  \right) \nonumber\\
		&=\Phi (\ft,k) \left( R_{1,k} - (\qmk{k} - R_{1,k}) \qmk{k}^{-1} R_{2,k} \right) .
	\end{align*}
	Thus, \eqref{eq:R1_R2} is equivalent to
	\begin{equation*}
		J_k := R_{1,k} + R_{1,k} \qmk{k}^{-1} R_{2,k} - R_{2,k} = 0 .
	\end{equation*}

	Hereafter, we show that
	\begin{equation}\label{eq:J_transition}
		J_{k+1} = A_k J_k A_{Q,k}^\top, \ k\in \bbra{0,\ft-1} .
	\end{equation}
    The above relationship with $ J_\ft = 0 $ and the invertibility of $ A_k $ and $ A_{Q,k} $ imply \eqref{eq:R1_R2}, which concludes \eqref{eq:final_coef}.
	For the left-hand side of \eqref{eq:J_transition}, it holds
	\begin{align}
		&J_{k+1} = A_k R_{1,k} A_k^\top - B_k B_k^\top + A_k R_{1,k} \qmk{k}^{-1} R_{2,k} A_{Q,k}^\top \nonumber\\
        &- A_k R_{1,k} \qmk{k}^{-1} A_k^{-1} B_k B_k^\top - B_k B_k^\top A_k^{-\top} \qmk{k}^{-1} R_{2,k} A_{Q,k}^\top \nonumber\\
		& + B_k B_k^\top A_k^{-\top} \qmk{k}^{-1} A_k^{-1} B_k B_k^\top - A_{Q,k} R_{2,k} A_{Q,k}^\top + B_{Q,k} B_{Q,k}^\top . \label{eq:J_next}
	\end{align}
	Here, we have
	\begin{align}
		B_k B_k^\top A_k^{-\top} \qmk{k}^{-1} R_{2,k} A_{Q,k}^\top &= B_{Q,k}B_{Q,k}^\top A_{Q,k}^{-\top} \qmk{k}^{-1} R_{2,k} A_{Q,k}^\top \nonumber\\
		&= B_{Q,k}B_{Q,k}^\top \qmk{k+1}^{-1} A_k R_{2,k} A_{Q,k}^\top \nonumber\\
		& = (A_k - A_{Q,k}) R_{2,k} A_{Q,k}^\top , \label{eq:J_next_5th}
	\end{align}
	where in the last line we used $ A_{Q,k} = A_k - B_{Q,k} B_{Q,k}^\top \qmk{k+1}^{-1} A_k $.
	On the other hand,
	\begin{align}
		&B_k B_k^\top A_k^{-\top} \qmk{k}^{-1} A_k^{-1} B_k B_k^\top \nonumber\\
        &= B_k B_k^\top (\qmk{k+1} + B_k B_k^\top)^{-1} B_k B_k^\top \nonumber\\
		&= - B_k (- I + I - B_k^\top (\qmk{k+1} + B_k B_k^\top)^{-1} B_k) B_k^\top \nonumber\\
		&= B_k B_k^\top - B_k (I + B_k^\top \qmk{k+1}^{-1} B_k)^{-1} B_k^\top \nonumber\\
		&= B_k B_k^\top - B_{Q,k} B_{Q,k}^\top . \label{eq:BB_BQBQ}
	\end{align}
	Moreover,
	\begin{align}
		&A_k R_{1,k} A_k^\top - A_k R_{1,k} \qmk{k}^{-1} A_k^{-1} B_k B_k^\top \nonumber\\
        &= A_k R_{1,k} A_k^\top - A_k R_{1,k} \qmk{k}^{-1} A_{Q,k}^{-1} B_{Q,k} B_{Q,k}^\top  \nonumber\\
		&= A_k R_{1,k} A_k^\top - A_k R_{1,k} A_k^\top \qmk{k+1}^{-1} B_{Q,k} B_{Q,k}^\top \nonumber\\
		&= A_k R_{1,k} A_{Q,k}^\top . \label{eq:J_next_1st_4th}
	\end{align}
	By \eqref{eq:J_next}, \eqref{eq:J_next_5th}, \eqref{eq:BB_BQBQ}, \eqref{eq:J_next_1st_4th}, it holds
	\begin{align*}
		J_{k+1} &= A_k R_{1,k} A_{Q,k}^\top + A_k R_{1,k} \qmk{k}^{-1} R_{2,k} A_{Q,k}^\top - A_k R_{2,k} A_{Q,k}^\top  \nonumber\\
		&= A_k J_k A_{Q,k}^\top .
	\end{align*}
	Hence, we come to \eqref{eq:final_coef}.

	By using \eqref{eq:final_coef}, we get \eqref{eq:bridge_A} as follows.
	\begin{align*}
		&\what{A}_{Q,k} = \bigl[ I - B_k B_k^\top A_k^{-\top} (\Phi (\ft,k) G_\rmc (\ft,k))^\dagger \nonumber\\
        &\qquad\qquad \times \Phi_\qm (\ft,k+1) \bigr] A_{Q,k} \nonumber\\
		&= A_k - B_{Q,k} B_{Q,k}^\top \qmk{k+1}^{-1} A_k - B_k B_k^\top A_k^{-\top} \nonumber\\
        &\quad \times \left(\Phi (\ft,k) G_\rmc (\ft,k)\right)^\dagger \Phi (\ft,k) \Phi (k,\ft) \qmk{\ft} \Phi (k,\ft)^\top \qmk{k}^{-1} \nonumber\\
		&= A_k - B_{Q,k} B_{Q,k}^\top \qmk{k+1}^{-1} A_k - B_k B_k^\top A_k^{-\top} \nonumber\\
        &\quad \times (\Phi (\ft,k) G_\rmc (\ft,k))^\dagger \Phi (\ft,k) (I - G_\rmc (\ft,k) \qmk{k}^{-1}) \nonumber\\
		&= A_k - B_{Q,k} B_{Q,k}^\top \qmk{k+1}^{-1} A_k - B_k B_k^\top A_k^{-\top} \nonumber\\
        &\quad \times (\Phi (\ft,k) G_\rmc (\ft,k))^\dagger \Phi (\ft,k) + B_k B_k^\top A_k^{-\top} \qmk{k}^{-1} \nonumber\\
		&= \what{A}_k .
	\end{align*}
	Lastly, by using \eqref{eq:final_coef} and \eqref{eq:BB_BQBQ}, we show \eqref{eq:bridge_cov} as
	\begin{align*}
		&\Lambda_{Q,k} = B_{Q,k} B_{Q,k}^\top - B_k B_k^\top A_k^{-\top} (\Phi (\ft,k) G_\rmc (\ft,k))^\dagger \nonumber\\
        &\qquad\quad \times  \Phi_\qm (\ft,k+1) B_{Q,k} B_{Q,k}^\top \nonumber\\
		&= B_k B_k^\top - B_k B_k^\top A_k^{-\top} \qmk{k}^{-1} A_k^{-1} B_k B_k^\top - B_k B_k^\top A_k^{-\top} \nonumber\\
        &\quad\times (\Phi (\ft,k) G_\rmc (\ft,k))^\dagger \Phi (\ft,k) (I - G_\rmc (\ft,k) \qmk{k}^{-1}) \nonumber\\
        &\quad \times A_k^{-1} B_k B_k^\top \nonumber\\
		&= B_k B_k^\top \nonumber\\
        &\quad - B_k B_k^\top A_k^{-\top} (\Phi(\ft,k) G_\rmc (\ft,k))^\dagger \Phi (\ft,k+1) B_k B_k^\top \nonumber\\
		&= \Lambda_k .
	\end{align*}
    In summary, we obtain (ii). The expression \eqref{eq:separate_KL} with (i),~(ii) completes the proof.

\bibliographystyle{IEEEtran}
\bibliography{TAC_maxent_bridge}

\begin{IEEEbiography}[{\includegraphics[width=1in,height=1.25in,clip,keepaspectratio]{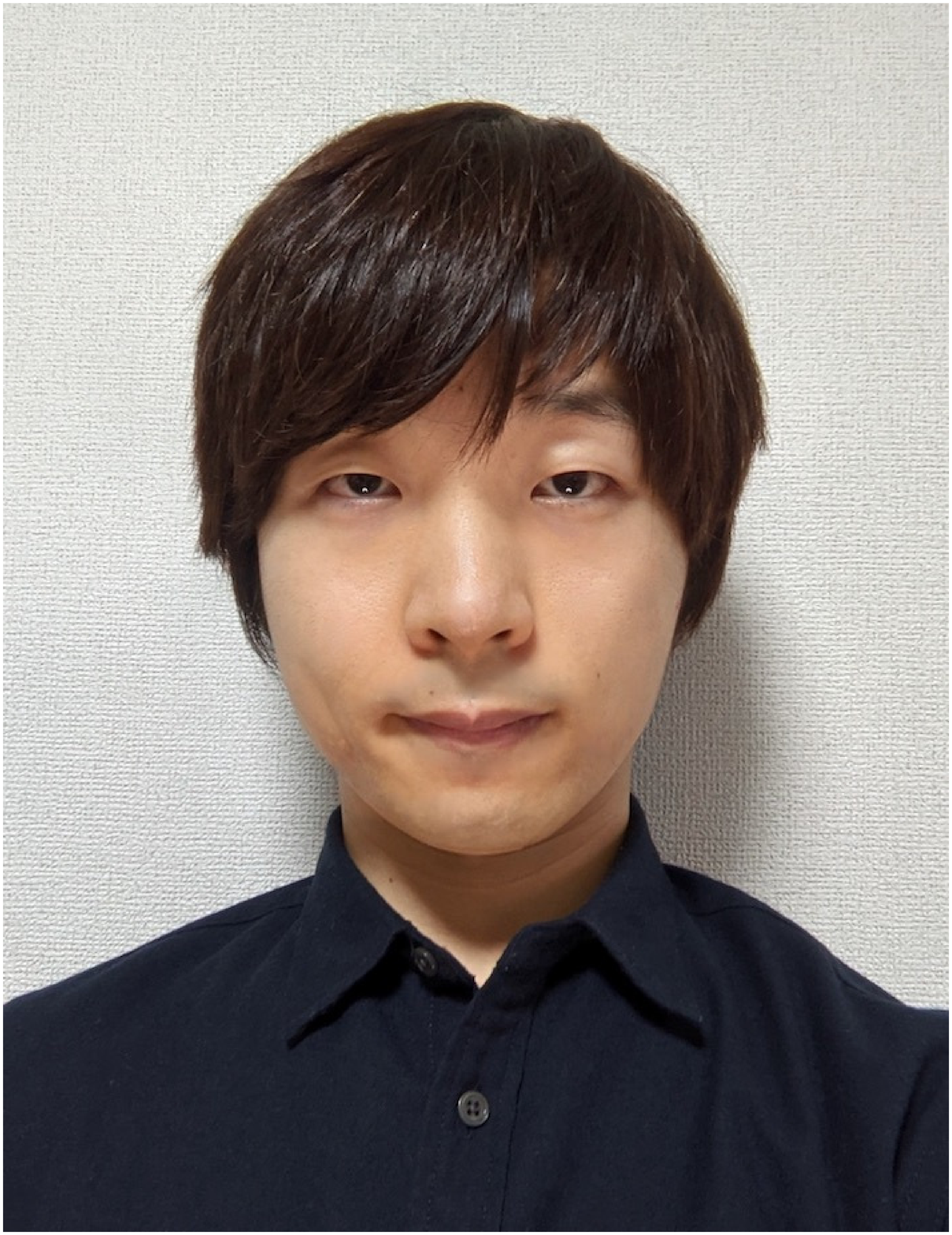}}]{Kaito Ito}
	(Member, IEEE) received his Bachelor's degree in Engineering and Master's degree and Doctoral degree in Informatics from Kyoto University in 2017, 2019, and 2022, respectively. 
	
	He is currently an Assistant Professor with the School of Computing, Tokyo Institute of Technology. His research interests include stochastic control, large-scale systems, and statistical learning.
\end{IEEEbiography}

\begin{IEEEbiography}[{\includegraphics[width=1in,height=1.25in,clip,keepaspectratio]{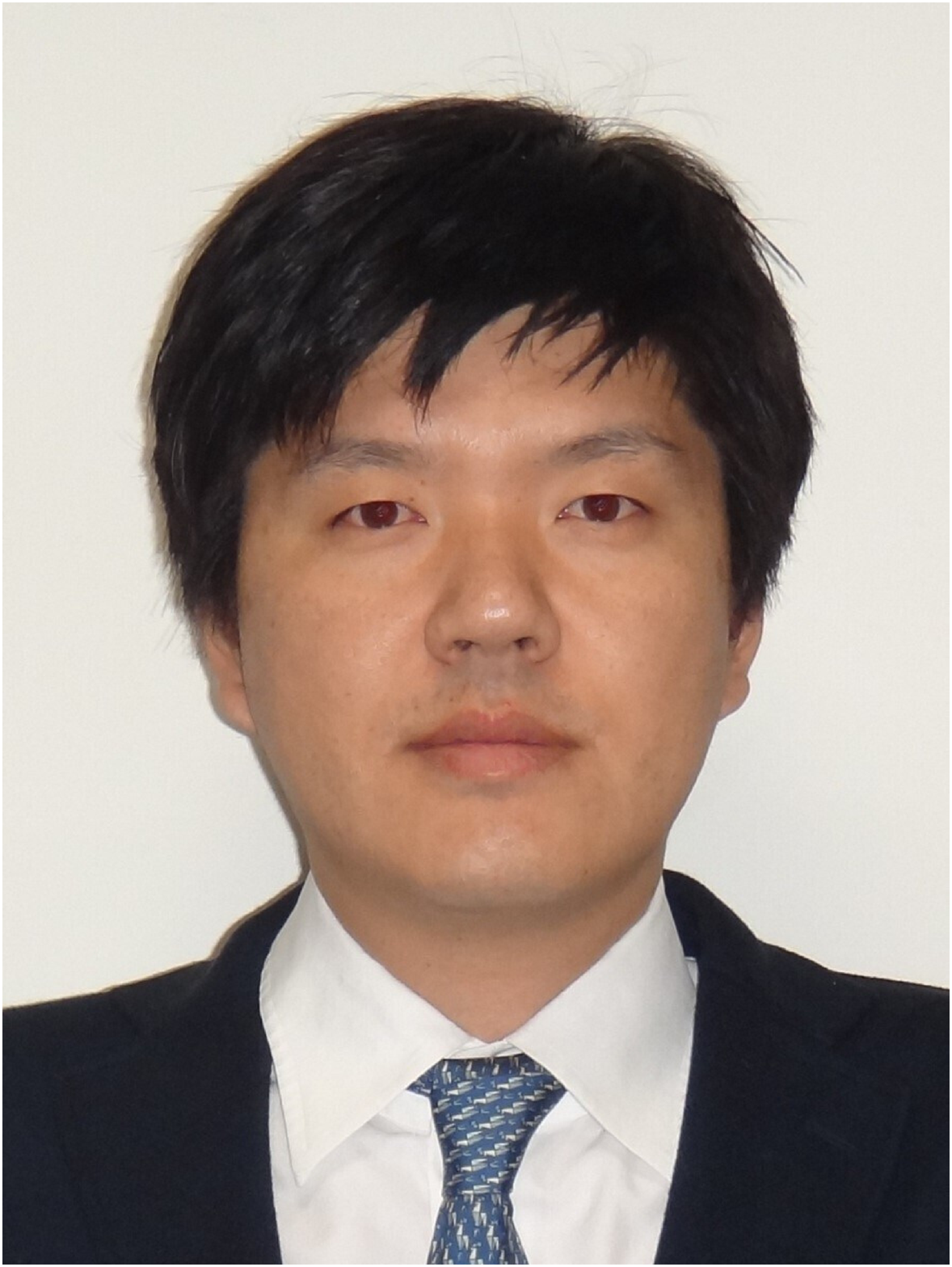}}]{Kenji Kashima}
	(Senior Member, IEEE) received his Doctoral degree in Informatics from Kyoto University in 2005. 
	
	He was with Tokyo Institute of Technology, Universit\"{a}t Stuttgart, Osaka University, before he joined Kyoto University in 2013, where he is currently an Associate Professor. His research interests include control and learning theory for complex (large scale, stochastic, networked) dynamical systems, as well as its applications. 
	He has served as an Associate Editor of IEEE Transactions of Automatic Control, IEEE Control Systems Letters, the IEEE CSS Conference Editorial Board and Asian Journal of Control.
\end{IEEEbiography}

\end{document}